\documentclass[12pt,a4paper,
 ]{amsart}

\usepackage{soul}
\usepackage{lscape}
\usepackage{graphicx}
\usepackage[all]{xy}

\usepackage{DLdef}

\makeatletter

\newcommand\surtab{\small\centering\extrarowheight=1\p@ \def\arraystretch{1.4}\tabcolsep3.6\p@ \doublerulesep=0.09\p@%
     \vskip 3\p@\advance\belowcaptionskip3\p@ \mathsurround\z@ \hyphenpenalty\z@ \doublehyphendemerits\z@ }
\newcommand\tcaption[2][0pt]{%
   \def\@tempa{#2}%
   \addtocounter{table}{1}\surtab\edef\@currentlabel{\thetable}%
   \ifx\@tempa\@empty \hfill \textbf{\tablename}~\thetable\kern#1 \else\textbf{\tablename}~\thetable.\kern 1ex{\boldmath\bfseries #2}\fi
}

\makeatother

\begin{document}

\graphicspath{{./S3/}{./li3/}}

\markboth{D. Chapovalov, M. Chapovalov, A. Lebedev, D. Leites}
{Almost affine Lie superalgebras}

\title
 {The classification of almost affine (hyperbolic) Lie superalgebras}

\author{Danil Chapovalov${}^1$, Maxim Chapovalov${}^2$,\\ Alexei Lebedev${}^3$,
Dimitry Leites${}^2$}

\address{${}^1$iBusiness, Stockholm, Sweden}

\address{${}^2$Dept. of Math., Stockholm University, Roslagsv. 101,
Kr\"aftriket hus 6, SE-106 91, Stockholm, Sweden;
mleites@math.su.se}

\address{${}^3$Nizhegorodskij University, pr. Gagarina 23, Nizhny
Novgorod, \\
RU-603950 Russia; yorool@mail.ru}

\thanks{DL is thankful to \'E. Vinberg and O. Shwartsman for helpful
comments and to A.~Protopopov for his help with our graphics, see
\cite{Pro}.}



\begin{abstract}
We say that an indecomposable Cartan matrix $A$ with entries in the
ground field of characteristic 0 is {\it almost affine} if the Lie
sub(super)algebra determined by it is not finite dimensional or
affine (Kac--Moody) but the Lie (super)algebra determined by any
submatrix of $A$, obtained by striking out any row and any column
intersecting on the main diagonal, is the sum of finite dimensional
or affine Lie (super)algebras. A Lie (super)algebra with Cartan
matrix  is said to be {\it almost affine} if it is not finite
dimensional or affine (Kac--Moody), and all of its Cartan matrices
are almost affine.

We list all almost affine Lie superalgebras over complex numbers
correcting two earlier claims of classification and make available
the list of almost affine Lie algebras obtained by Li Wang Lai.
\end{abstract}

\maketitle

\section{Introduction}

In what follows the ground field is $\Cee$; all Cartan matrices $A$
are supposed to be indecomposable and {\bf normalized} so that
$A_{ij}= 0$ or 1 or 2 (for details, see eq. \eqref{norm}); besides
$A_{ij}= 0\Longleftrightarrow A_{ji}= 0$;  for a given $n\times n$
matrix, let its {\it size} be $n$.

We say that a given Cartan matrix $A$ with entries in the ground
field is {\it almost affine} if the Lie sub(super)algebra determined
by it is not finite dimensional or affine (Kac--Moody) but the Lie
(super)algebra determined by any submatrix of $A$, obtained by
striking out any row and any column intersecting on the main
diagonal, is the sum of finite dimensional or affine Lie
(super)algebras. A Lie (super)algebra with Cartan matrix  is said to
be {\it almost affine} if it is not finite dimensional or affine
(Kac--Moody), and all of its Cartan matrices are almost affine.

{\bf Warning}: A given Lie superalgebra which is not almost affine
may possess an almost affine Cartan matrix which goes into a not
almost affine matrix under an odd reflection. We are interested in
the properties of Lie superalgebra, not in those of one or several
of its Cartan matrices.

\ssec{Motivations. Our result} Under the name \lq\lq hyperbolic"
and \lq\lq overextended" the almost affine Lie algebras find
applications in \lq\lq a variety of physical models in
two-dimensional field theories (supergravity, string theory,
cosmological billiards) \cite{GN, GNW, HM, H, DHN, DN, KN1, KN2,
BS}" (quotation from \cite{FS} with references updated but without
any attempt to give an exhaustive review of the literature).

There is another type of Lie (super)algebras that used to go under
the same name \lq\lq hyperbolic", but are defined differently and
currently are referred to as {\it Lorentzian} Lie algebras; for
their precise definition, see \cite{RU, GN}).

The study of Lorentzian Lie algebras makes superization not just
natural, but rather \so{inevitable} we'd say, see  \cite{RU, GN}:
Borcherds, and later Gritsenko and Nikulin found various
applications of simple Lorentzian  Lie algebras and {\bf
super}algebras (for one of these applications Borcherds was awarded
with a Fields medal). This certainly justifies the quest for the
simple Lorentzian  Lie algebras and {\bf super}algebras.

Whereas the almost affine Lie algebras are useful, e.g., for
cosmologic billiards, at the moment the only applications of almost
affine Lie {\bf super}algebras we know of are due to the fact that
some of them (those of rank 3) coincide with the known Lorentzian
Lie {\bf super}algebras. This is already good, but the notion of
almost affine Lie (super)algebras seems to be most natural even
without such coincidence, and hence worth investigating. Besides,
any almost affine Lie superalgebra whose even part is an almost
affine Lie algebra might hint at a hidden supersymmetry of the
problem related with the latter.

Although the classification problem of the almost affine Lie {\bf
super}algebras (twice claimed to be done) should be (and is) much
simpler than the classification problem of Lorentzian Lie
(super)algebras (still an open problem, and it is not even clear if
it is tame even if there are finitely many of them), it was not
solved in one go.

This paper was written after we failed to understand even the
definitions given in \cite{FS} to say nothing of results (which, in
turn, were supposed to correct the results of \cite{TDP}); several
counterexamples to the claims of \cite{FS} immediately spring to
mind. Here we rectify the results of \cite{TDP} and \cite{FS}; in
particular, we give precise definitions and an algorithm to provide
with the complete classification of almost affine Lie superalgebras.

Our result is obtained with the help of a computer and
double-checked manually. We do not list an intermediate result ---
classification of almost affine Cartan matrices in the super case.

\ssec{On linguistics} {\it A posteriori} it turns out that almost
affine Lie algebras (nothing super) with symmetrizable and integer
Cartan matrix of size $>2$ are what was lately called {\it
hyperbolic} Lie algebras because their Weyl groups are hyperbolic
Coxeter groups. In the absence of an {\it adequate} super version of
the notion of the Weyl group, the term \lq\lq hyperbolic" is
meaningless in super setting.

There is another unfortunate term --- \lq\lq overextended", meaning,
actually, \lq\lq twice extended" or \lq\lq doubly extended". This
means that, having extended the Dynkin {\bf graph}, we extend it
still further to get the Dynkin graph of an almost affine Lie
algebra. However, the term \lq\lq extended Lie algebra" is occupied
and means something different.

Actually, even the definition of \lq\lq hyperbolic"  Lie algebras
(nothing super) by means of Cartan matrices with integer
coefficients does not look natural; Borcherds was the first to rebel
against the fossilized definition (with a remarkable success
tempting one to go on).

\ssec{On setting of the problem} Simple finite dimensional Lie
algebras need not be introduced:  during the past century, they
proved their usefulness in mathematics and physics on numerous
occasions. Each of these algebras possesses a number of useful
properties (has a symmetrizable Cartan matrix, and hence an
invariant symmetric bilinear form; invertibility of its Cartan
matrix implies non-degeneracy of this form; possesses various
gradings, and so on).

Occasional applications of infinite dimensional Lie algebras that
possess some of the above properties (and ensuing  applications)
motivated a quest for classifications retaining some of the
properties (since it is impossible to retain all of them); these
classifications resulted in rich theories of at least two of the
following three types of infinite dimensional Lie algebras; even
more of new applications were found. The three types of algebras we
have in mind are (for applications, see references in parentheses)
\begin{enumerate}\label{1-3}
  \item Lie algebras of vector fields whose coefficients are polynomials,
formal series, or Laurent polynomials (\cite{Gr, LSa}; for
superization, see \cite{GLS1, LSh, Kac});
  \item affine Kac-Moody algebras (\cite{K}; for superization of the list, see
  \cite{FLS, GL1});
  \item analogs of the Lie algebra of \lq\lq matrices of complex
  size" (see \cite{GL2} and references therein).
\end{enumerate}

The only  catch was (and still is): what are the criteria for \lq\lq
similarity" of properties?

Which properties of finite dimensional simple Lie algebras should we
(try to) retain and to what extent? What are the new properties of
certain infinite dimensional Lie algebras that are \lq\lq useful" in
applications and that finite dimensional simple Lie algebras do not
possess?

Already the above examples 1--3 manifestly show that there are
\so{several} classes of reasonably \lq\lq nice" infinite dimensional
Lie algebras \lq\lq similar", to an extent, to simple finite
dimensional Lie algebras. And it is immediately clear that the
property very useful for classification purposes  --- simplicity ---
is not perfect in applications: certain \lq\lq relatives" of simple
Lie algebras (such as their non-trivial central extensions,
deformations, certain derivation algebras) are more useful.

Lie {\bf super}algebras, which mathematicians started to consider in
1930s, became a topic of conscious interest and deliberate study in
mid-1970s and all the above applies to them as well.

In the above discussion we assumed that the ground field is $\Cee$
(although the papers studying infinite dimensional Lie
(super)algebras over fields of prime characteristic already started
to appear, this area is not ready to be studied yet, we think).
However, in physical applications, real forms of complex algebras
are often more natural objects; classification of \lq\lq nice" Lie
(super)algebras over $\Ree$ is a natural and important ramification
of the classification problem over $\Cee$. Such classifications are
not as obvious as in the case of finite dimensional simple Lie
algebras over $\Cee$ (cf. \cite{FH} with, e.g., a difficult result
of \cite{CK} on classification of real forms of simple Lie algebras
of polynomial vector fields classified over $\Cee$ by Leites and
Shchepochkina \cite{LSh}, or with unexpected Serganova's results on
classification of real forms of affine Kac--Moody Lie
(super)algebras \cite{FLS} and stringy Lie (super)algebras
\cite{Se2}).

{\bf Superization} of everything gives one more dimension to the
classification problem. There are usually two types of super
notions: a straightforward one and more involved ones. The most
interesting are, of course, the involved ones, especially if they
naturally appear. This is precisely the case in the problem we are
going to consider!

We will carefully distinguish between the distinct types of algebras
(and insist on different names: {\it Lorentzian} (considered in
\cite{RU, GN}) and {\it almost affine}; we also suggest to never use
the adjective {\it hyperbolic} speaking about Lie (super)algebras
while keeping it for Coxeter groups) but first recall where these
terms originate from.

Every simple finite dimensional Lie algebra is of the form $\fg(A)$
(where $A$ is called a {\it Cartan matrix}, for the definition, see
\S2); moreover
\begin{enumerate}
  \item the Cartan matrix $A$ of $\fg(A)$ is {\it integer} (after a
normalization to be described; hereafter all Cartan matrices are
normalized unless otherwise stated);
  \item the Cartan matrix $A$ of $\fg(A)$ is
 {\it invertible};
  \item $\fg(A)$ is $\Zee$-graded, and even $\Zee^n$-graded, where
  $n$ is the size of $A$.
\end{enumerate}

There are also other characterizations of these matrices\footnote{Or
Lie algebras $\fg(A)$, which at this stage seems to be the same
thing but will turn to be something else for almost affine {\bf
super}algebras.} $A$, perhaps, no less interesting.

\ssec{Linguistics continued} 
A submatrix $B$ of $A$ obtained by striking out a rows
and a column with same numbers is said to be {\it principal}.

\sssbegin{Lemma}[\cite{Vi}] Let $A$ be a real matrix (normalized so
that $A_{ii}=2$ or $0$ for all $i$) such that

1) $A$ is indecomposable (i.e., can not be reduced to a
block-diagonal form by a permutation of its rows and same columns),

2) $A_{ij}\leq 0$ for $i\neq j$,

3) $A_{ij}= 0\Longleftrightarrow A_{ji}= 0$.

\noindent Then $A$ satisfies one and only one of the following
conditions:

Fin) $\det B>0$ for every principal submatrix $B$ of $A$, and $\det
A>0$;

Aff) $\det A =0$ and $\det B>0$ for every principal submatrix $B$ of
$A$;

Ind) otherwise.
\end{Lemma}

\sssec{Comments to definitions} The Lie algebra $\fg(A)$ (and its
Cartan matrix $A$) is said to belong to the {\it finite}, {\it
affine} and {\it indefinite type}, in accordance with the cases of
Lemma. In particular, although we are dealing with complex algebras
the above Lemma requires that the entries of $A$ should be real (to
make inequalities meaningful). There is no apparent reason to impose
this restriction from the very beginning (except that over $\Cee$ we
do not know what might the analog of the condition 2) of Lemma and
its reformulation be), moreover, in the super case, there is no
analog of this Lemma, anyway.

In cases Fin and Aff, the matrices $A$ are symmetrizable and the
off-diagonal entries are non-positive integers.

Kac and, independently, Moody singled out a class of Cartan matrices
$A$ such that the Lie algebras $\fg(A)$ recovered from $A$ by
(almost) the same rules as for finite dimensional Lie algebras
possess many of the properties of their finite dimensional models.
Kac  later developed a theory (for a summary, see \cite{K}) of these
Lie algebras nowadays called {\it affine Kac-Moody algebras}.

For the shearing parameter convenient to single out the \lq\lq nice"
objects among $\Zee$-graded Lie algebras
$\fg=\mathop{\oplus}\limits_{i\in\Zee}\fg_i$ such that $\dim
\fg_i<\infty$ for all $i$, Kac selected the {\it growth}
\begin{equation}\label{grt}
    \text{grt}(\fg):=\mathop{\limsup}\limits_{n \tto \infty}
    \frac{\ln \dim\mathop{\oplus}\limits_{|i|\leq n}\fg_i}{\ln n}.
\end{equation}
Kac proved that growth gives a simple and precise separation:
\begin{equation}\label{grt1}
\text{grt}(\fg)=\begin{cases}0&\text{if $\dim\fg<\infty
$};\\
1&\text{if $\fg$ is an affine Kac-Moody algebra};\\
\infty&\text{otherwise}.\end{cases}
\end{equation}

Soon, however, it became clear that, among the sea of Lie algebras
of the form $\fg(A)$ of infinite growth and apparently unyielding to
study, there are certain islands with nice properties.

\sssec{Hyperbolic Coxeter groups and almost affine Lie
algebras}\label{sshyp} According to Bourbaki (\cite{Bo}) a Coxeter
group is said to be of a {\it (compact) hyperbolic} type if any its
subgroup generated by any proper subset of generators is an (affine)
Weyl group. These groups act on the hyperbolic (or Lorentzian)
space; they are compact if so is the fundamental domain. In this
definition, both the adjectives \lq\lq hyperbolic" and \lq\lq
compact" are meaningful: the other types being elliptic (or
spherical) corresponding to the Coxeter groups acting on spheres,
and parabolic (or Euclidean) corresponding to the Coxeter groups
acting on $\Ree^n$, 
see \cite{Bo}.

The Bourbaki adjectives for Coxeter groups were (thoughtlessly, we
think) extended to the Lie algebras with these groups as their Weyl
groups. In the process, the reasonable adjective \lq\lq compact" was
replaced by meaningless \lq\lq strict". The adjective \lq\lq
hyperbolic" still retains (no doubt coincidentally) some sense for
Lie algebras of rank $>2$: if the Cartan matrix of an almost affine
Lie algebra is symmetrizable, than it can be considered as the Gram
matrix of the non-degenerate bilinear form of Lorentzian signature.

Moody \cite{Moo} uses the adjective \lq\lq hyperbolic" speaking
about root systems, which is reasonable but --- and this is vital
for us --- unclear how to superize (except, to an extent, revealed
by Serganova, \cite{Se3}). But we intend to deal with Lie
superalgebras, not their root systems, anyway.

Another approach to superization of \lq\lq hyperbolic" Lie algebras
is by looking at their Weyl groups. But Lie superalgebras have no
Weyl groups (except in a few cases)! So, although for Lie
superalgebras, one can define {\it neighboring systems of simple
roots} and reflections that interchange these neighboring systems,
we prefer to deal with notions better understood than super analogs
of Weyl groups, namely, with Lie superalgebras themselves.

The classification of almost affine Lie algebra $\fg(A)$ shows that
\begin{equation}\label{int}
\text{all off-diagonal entries of $A$ are non-positive integers if
size$(A)>2$}.
\end{equation}
It is clear, however, that to require \eqref{int} {\it a priori} is
unjustified; it is more natural to consider the entries of $A$
belonging to the ground field. Borcherds was, perhaps, the first to
demonstrate usefulness of slackening the requirement \eqref{int},
cf. the review \cite{GN}, the book \cite{RU} and references therein.

The definition of almost affine Lie algebras immediately implies
their description for indecomposable Cartan matrices of size 2: all
matrices of the following form will do (the excluded cases being
finite dimensional or affine):
\begin{equation}\label{hysize2}
\begin{lmatrix}2&a\\b&2\end{lmatrix},\text{~~where $a,b\in\Cee\setminus{0}$,
$\begin{matrix}\\
(a,b)\neq &(-1, -2), (-2, -1),(-1, -3), (-3, -1),\\&(-1, -4),(-4,
-1), (-1,-1), (-2, -2).\end{matrix}$}
\end{equation}
Bourbaki  \cite{Bo} (Ex. 26 to \S1 of Ch.4) describes simple almost
affine Lie algebras of rank $>2$ (the implicit result is due to
L.~Solomon \cite{So} and R.~Steinberg \cite{Ste}, 1966): There are
finitely many of them, in particular, their rank is bounded: it is
$\leq 10$. In \cite{K}, this exercise from \cite{Bo} is repeated
with the answer still implicit although by the time the first
edition of the book appeared it was known:

In 1983, Kobayashi and Morita classified the almost affine Lie
algebras with indecomposable symmetrizable Cartan matrix of size
$>2$ \cite{KoMo}. In 1988, Li Wang Lai \cite{Li} obtained the {\it
complete} answer for Dynkin diagrams on $>2$ vertices: There are 238
almost affine Lie algebras; 142 of these algebras have a
symmetrizable Cartan matrix. Both \cite{Li} and \cite{KoMo} were
published in unpopular journals and did not put their results in
arXiv later, so one usually cites (sometimes trustingly, sometimes
to correct, see \cite{BS}) a later result by
C.~Sa{\c{c}}lio{\u{g}}lu \cite{S} who tackled the same problem.
Though unjust to \cite{Li} and \cite{KoMo}, this practice is
understandable: Sa{\c{c}}lio{\u{g}}lu's paper has gaps even in his
partial (concerning only symmetrizable Cartan matrices) answer but
it is very interestingly written. Besides, Chinese put their last
names in front of their first ones, which hampers proper citing, cf.
\cite{FS}. Since the results of Li are difficult to access but are
of interest, we reproduce them in Appendix.

\parbegin{Remark}\label{tafrek} To require  that $a, b\in-\Zee_+$ in
\eqref{hysize2} are positive integers such that $ab>4$ does not seem
natural and the result of \cite{KM}, where $b=a$ and which
analytically depends on $a$, supports our desire to let off-diagonal
elements belong to the ground field.
\end{Remark}

\ssec{Superization} Classification of simple finite dimensional Lie
superalgebras over $\Cee$ is a result of disjoint work of several
teams of researchers, see \cite{Kapp}. Of these, Kac was the first
to try to classify simple finite dimensional Lie superalgebras of
the form $\fg(A)$ over $\Cee$. Kac pointed at an important fact: For
certain simple Lie superalgebras (first discovered by Kaplansky
\cite{Kapp}),
\begin{equation}\label{param}
\text{the entries of $A$ can belong to the ground field even if
$\text{grt}(\fg(A))=0$}.
\end{equation}

The classification of simple Lie superalgebras of the form $\fg(A)$
of finite growth was almost correctly conjectured in \cite{FLS}, for
a missing item, see \cite{GL1}; the conjecture was proved for
symmetrizable matrices $A$ in \cite{FLS} (based on \cite{Se}) and
independently in \cite{vdL}; for non-symmetrizable matrices $A$, the
proof was only recently published \cite{HS}. Finally:
\begin{equation}\label{non-sy}
\text{\begin{minipage}[c]{14cm} If a simple Lie superalgebra is of
finite dimension, its Cartan matrix $A$ is symmetrizable, whereas,
unlike non-super case, if $\text{grt}(\fg(A))=1$, then $A$ can be
non-symmetrizable.
\end{minipage}}
\end{equation}

\section{What $\fg(A)$ is}\label{Sg(A)}

The finite dimensional Lie superalgebras of type $\fsl$ and $\fosp$
consist of linear operators preserving the volume element and the
even symmetric (and anti-symmetric) form, respectively. Some of them
are simple, some are not. For an interpretation of the exceptional
Lie superalgebras $\fa\fg(2)$ and $\fa\fb(3)$, see \cite{CE2}; but
nothing is more convenient (at least, for computers) than their
presentation in terms of Cartan matrices \cite{GL1}. The Lie
superalgebra $\fosp(4|2)$ admits a parametric family of
deformations, each possessing Cartan matrix, and the elements
$\fosp_\alpha(4|2)$ (where $\alpha\neq 0, 1$) of this family will be
sometimes denoted $\fd(\alpha)$ for brevity. For an interpretation
of the parameter $\alpha$, see \cite{BGL}. The following Lie
superalgebras are needed for the proof, but not the answer:

The Lie superalgebra $\fq(n)$ is a \lq\lq queer" analog of $\fgl$;
it preserves an odd complex structure given by an odd operator. The
indigenous \lq\lq queer" trace on it singles out the queertraceless
subalgebra $\fs\fq(n)$. The projectivization (passage to the
quotient modulo center consisting of scalar matrices) leads to
$\fp\fsl(n|n)$ and $\fp\fs\fq(n)$.

The Lie algebra $\fsvect_\alpha(1|2)$ consists of vector fields with
Laurent polynomials as coefficients and preserving the volume
element $t^\alpha\vvol(t|\xi_1,\xi_2)$.

\ssec{Warning: certain of $\fsl$'s and all $\fpsl$'s have no Cartan
matrix. Their relatives that have Cartan matrices}\label{warn} For
the most reasonable definition of Lie algebra with Cartan matrix
over $\Cee$, see \cite{K}. The same definition applies, practically
literally, to Lie superalgebras. However, the usual sloppy practice
is to attribute Cartan matrices to many of those (usually simple)
Lie superalgebras which, strictly speaking, have no Cartan matrix!

Although it may look strange for the reader with non-super
experience over $\Cee$, neither the simple Lie superalgebra
$\fpsl(a|a)$ nor the Lie superalgebra $\fsl(a|a)$ have Cartan
matrix.

Their relative possessing a Cartan matrix is $\fgl(a|a)$, and for
the grading operator we take $E_{1,1}$.

Since often all the Lie (super)algebras involved (the simple one,
its central extension, the derivation algebras thereof) are needed
(and only representatives of one of the latter types of Lie
(super)algebras are of the form $\fg(A)$), it is important to have
(preferably short and easy to remember) notation for each of them.
For example: $\fpsl$, $\fsl$, $\fp\fgl$ and $\fgl$.

\ssec{Generalities} Let us start with the construction of a Cartan
matrix Lie (super)al\-geb\-ra (in what follows briefly called CM Lie
(super)algebra or even CMLA or CMLSA). Let $A=(A_{ij})$ be an
$n\times n$-matrix. Let $\rk A=n-l$. It means that there exists an
$l\times n$-matrix $T=(T_{ij})$ such that
\begin{equation}\label{rankCM}
\begin{tabular}{l}
a) the rows of $T$ are linearly independent;\\
b) $TA=0$ (or, more precisely, \lq\lq zero $l\times n$-matrix'').
\end{tabular}
\end{equation}
Indeed, if $\rk A^T=\rk A=n-l$, then there exist $l$ linearly
independent vectors $x_i$ such that $A^Tx_i=0$; set
$$
T_{ij}=(x_i)_j.
$$

Let the elements $e_i^\pm$ and $h_i$, where $i=1,\dots,n$, generate
a Lie superalgebra denoted $\tilde \fg(A, I)$,\index{$\tilde\fg(A,
I):=\tilde \fg(A, I)$} where $I=(i_1, \dots i_n)\in(\Zee/2)^n$ is a
collection of parities ($p(e_j^\pm)=i_j$), free except for the
relations (here $\fh:=\Span(h_i)_{i=1}^n$)
\begin{equation}\label{gArel_0}
{}[e_{i}^+, e_{j}^-] = \delta_{ij}h_i; \quad [h_i, e_{j}^\pm]=\pm
A_{ij}e_{j}^\pm\text{~~for any  $i, j$};\quad \;[\fh, \fh]=0.
\end{equation}
The Lie (super)algebras with Cartan matrix that we denote by $
\fg(A, I)$\index{$\fg(A, I)$} are quotients of $\tilde \fg(A, I)$
modulo the ideal we explicitly describe in \cite{LCh}. Implicitly,
the additional to \eqref{gArel_0} are relations $R_i=0$ whose left
sides are implicitly described as
\begin{equation}
\label{myst} \text{\begin{minipage}[c]{14cm} \lq\lq the $R_i$ that
generate the maximal ideal $\fr$ among the ideals of $\tilde\fg (A,
I)$ whose intersection with the span of the above $h_i$ and the
$d_j$ described in eq. \eqref{central3}  is zero.
 \rq\rq\end{minipage}}
\end{equation}

Set
\begin{equation}\label{central}
c_i=\sum_{j=1}^n T_{ij}h_j, \text{~~where~~} i=1,\dots,l.
\end{equation}
Then, from the properties of the matrix $T$, we deduce that
\begin{equation}\label{central1}
\begin{tabular}{l}
a) the elements $c_i$ are linearly independent;\\
b) the elements $c_i$ are central, because\\
$[c_i,e_j^\pm]=\pm\left(\sum\limits_{k=1}^n T_{ik}A_{kj}\right)
e_j^\pm=\pm (TA)_{ij} e_j^\pm $.
\end{tabular}
\end{equation}

The existence of central elements means that the linear span of all
the roots is only $(n-l)$-dimensional. (This can be explained even
without central elements: The weights can be considered as
column-vectors with $i$-th element being the corresponding
eigenvalue of $\ad_{h_i}$. The weight of $e_i$ is the $i$-th column
of $A$. Since $\rk A=n-l$, the linear span of all columns of $A$ is
$(n-l)$-dimensional just by definition of the rank. Since any root
is an (integer) linear combination of the weights of the $e_i$, the
linear span of all roots is $(n-l)$-dimensional.) This means that
some elements which we would like to see having different weights
have, actually, identical weights. To fix this, we do the following:
Let $B$ be an arbitrary $l\times n$-matrix such that
\begin{equation}\label{matrixB}
\text{the~~}(n+l)\times
n\text{-matrix~~}\begin{lmatrix}A\\B\end{lmatrix}\text{~~has
rank~}n.
\end{equation}
Let us add to the algebra the grading elements $d_i$, where
$i=1,\dots,l$, subject to the following relations:
\begin{equation}\label{central3}
{}[d_i,e_j^\pm]=\pm B_{ij}e_j;\quad [d_i,d_j]=0;\quad [d_i,h_j]=0.
\end{equation}
The last two relations mean that the $d_i$ lie in the Cartan
subalgebra, and even in the maximal torus $\fh$; from now on we
assume that $\fh$ is spanned not only by the $h_i$ but also by the
$d_j$ for all $i$ and $j$.

Note that these $d_i$ are {\it outer} derivations of $\fg^{(1)}$,
i.e., they can not be obtained as linear combinations of brackets of
the elements of $\fg^{(1)}$.

\ssec{Roots and weights}\label{roots} In this subsection, $\fg$
denotes one of the algebras $\fg(A,I)$ or $\tilde{\fg}(A,I)$.

Let $\fh$ be the span of the $h_i$ and the $d_j$. The elements of
$\fh^*$ are called {\it weights}.\index{weight} For a given weight
$\alpha$, its {\it weight subspace} of $\fg$ is defined as
$$
\fg_\alpha=\{x\in\fg\mid [h, x]=\alpha(h)x\text{~ for all
$h\in\fh$}\}.
$$

Any non-zero element $x\in\fg_\alpha$ is said to be {\it of weight
$\alpha$}.

\sssbegin{Statement}[\cite{K}] The space $\fg$ can be represented as
a direct sum of subspaces
$$
\fg=\mathop{\bigoplus}\limits_{\alpha\in \fh^*} \fg_\alpha.
$$
\end{Statement}

By construction, the elements\footnote{By abuse of notation we
retain the notations $e_i^\pm$ and $h_i$ --- the elements of
$\tilde{\fg}(A,I)$ --- for their images in $\fg(A,I)$ and
$\fg^{(1)}(A,I)$.} $e_i^\pm$ with the same superscript (either $+$
or $-$) have linearly independent weights $\alpha_i$, and any
$\alpha$ such that $\fg_\alpha\neq 0$ lies in the $\Zee$-span of
$\{\alpha_1,\dots,\alpha_n\}$.

The algebra $\fg$ has also an $\Ree^n$-grading such that $e_i^\pm$
has grade $(0,\dots,0,\pm 1,0,\dots,0)$, where $\pm 1$ stands in the
$i$-th slot (this can also be considered as $\Zee^n$-grading, but we
use $\Ree^n$ for simplicity of formulations). This grading is
equivalent to the weight grading of $\fg$. We will identify the
degrees with respect to these two gradings (this is used in
\eqref{x}).

Any non-zero element $\alpha\in\Ree^n$ is called {\it a
root}\index{root} if the corresponding eigenspace of grade $\alpha$
(which we denote $\fg_\alpha$ by abuse of notation) is non-zero. The
set $R$ of all roots is called {\it the root system}\index{root
system} of $\fg$.

Clearly, the subspaces $\fg_\alpha$ are purely even or purely odd,
and the corresponding roots are said to be \textit{even} or
\textit{odd}.

\ssec{Systems of simple and positive roots} In this subsection,
$\fg=\fg(A,I)$, and $R$ is the root system of $\fg$.

For any subset $B=\{\sigma_{1}, \dots, \sigma_{m}\} \subset R$, we
set (we denote by $\Zee_{+}$ the set of non-negative integers):
$$
R_{B}^{\pm} =\{ \alpha \in R \mid \alpha = \pm \sum n_{i}
\sigma_{i},\;\;n_{i} \in \Zee_{+} \}.
$$
The set $B$ is called a {\it system of simple roots~}\index{Root!
simple system of}\; of $R$ (or $\fg$) if $ \sigma_{1}, \dots ,
\sigma_{m}$ are linearly independent and $R=R_B^+\cup R_B^-$. Note
that $R$ contains basis coordinate vectors, and therefore spans
$\Ree^n$; thus, any system of simple roots contains exactly $n$
elements.

Let $(\cdot,\cdot)$ be the standard Euclidean inner product in
$\Ree^n$. A subset $R^+\subset R$ is called a {\it system of
positive roots~}\index{Root! positive system of}\; of $R$ (or $\fg$)
if there exists $x\in\Ree^n$ such that
\begin{equation}\label{x}
\begin{split}
 &(\alpha,x)\in\Ree\backslash \{0\}\text{ for all $\alpha\in R$},\\
 &R^+=\{\alpha\in R\mid (\alpha,x)>0\}.
\end{split}
\end{equation} Since $R$ is a finite (or infinite but countable) set, then the set
$$\{y\in\Ree^n\mid\text{there exists $\alpha\in R$ such that
} (\alpha,y)=0\}
$$ is a finite/countable union of $(n-1)$-dimensional
subspaces in $\Ree^n$, so it has zero measure. So for almost every
$x$, condition (\ref{x}) holds.

By construction, any system $B$ of simple roots is contained in
exactly one system of positive roots, which is precisely $R_B^+$.

\sssbegin{Statement} Any finite system $R^+$ of positive roots of
$\fg$ contains exactly one system of simple roots. This system
consists of all the positive roots (i.e., elements of $R^+$) that
can not be represented as a sum of two positive
roots.\end{Statement}

We can not give an {\it a priori} proof of the fact that each set of
all positive roots each of which is not a sum of two other positive
roots consists of linearly independent elements. This is, however,
true for finite dimensional Lie algebras and superalgebras of the
form $\fg(A)$.

\ssec{Normalization convention}\label{normA} Different (equivalent)
pairs $(A_B,I_B)$ may correspond to a given system of simple roots.
It would be nice to find a convenient way to fix some distinguished
pair $(A_B,I_B)$ in the equivalence class. For the role of the
\lq\lq best'' (first among equals) order of indices we propose the
one that minimizes the value
\begin{equation}\label{minCM}
\max\limits_{i,j\in\{1,\dots,n\}\text{~such that~}(A_B)_{ij}\neq
0}|i-j|,
\end{equation}
i.e., gather the non-zero entries of $A$ as close to the main
diagonal as possible. Observe that for the Lie algebras of type $E$
the standard (Bourbaki) numbering differs from \eqref{minCM}.

Clearly,
\begin{equation}
\label{rescale} \text{the rescaling
$e_i^\pm\mapsto\sqrt{\lambda_i}e_i^\pm$, sends $A$ to $A':=
\diag(\lambda_1, \dots , \lambda_n)\cdot A$.} 
\end{equation}

Two pairs $(A, I)$ and $(A', I')$ are said to be {\it equivalent} if
$(A', I')$ is obtained from $(A, I)$ by a composition of a
permutation of parities and a rescaling $A' = \diag (\lambda_{1},
\dots, \lambda_{n})\cdot A$, where $\lambda_{1}\dots \lambda_{n}\neq
0$. Clearly, equivalent pairs determine isomorphic Lie
superalgebras.

The rescaling affects only the matrix $A_B$, not the set of parities
$I_B$. The Cartan matrix $A$ is said to be {\it
normalized}\index{Cartan matrix, normalized} if
\begin{equation}
\label{norm} A_{jj}=0\text{~~ or 1, or 2 (only if $i_j=\ev$).}
\end{equation}
In order to distinguish the cases $i_j=\ev$ from $i_j=\od$, we write
$A_{jj}=\ev$ or $\od$, instead of 0 or 1, if $i_j=\ev$. \so{We will
only consider normalized Cartan matrices; for them, we do not have
to describe $I$.}

The row with a 0 or $\ev$ on the main diagonal can be multiplied by
any nonzero factor; we usually multiply it so as to make $A_{B}$
symmetric, if possible.

\ssec{Equivalent systems of simple roots} \label{EqSSR} Let
$B=\{\sigma_1,\dots,\sigma_n\}$ be a system of simple roots. Choose
non-zero elements $\tilde{e}_i^\pm\in\fg_{\pm\sigma_i}$; set $\tilde
h_{i}=[ \tilde{e}_{i}^{+}, \tilde{e}_{i}^-]$ , $A_{B} =(A_{ij})$,
where $A_{ij} =\sigma_{i}(\tilde{h_{j}})$ and
$I_{B}=\{p(\tilde{e}_{1}), \cdots, p(\tilde{e}_{n})\}$. (The pair
$(A_B,I_B)$ constructed here is not uniquely defined by $B$, but all
the pairs $(A_B,I_B)$ are equivalent to each other, and for any such
pair $(A_B,I_B)$, we have $\fg(A_B, I_B)\simeq\fg(A,I)$.) 

Two systems of simple roots $B_{1}$ and $B_{2}$ are said to be {\it
equivalent} if the pairs $(A_{B_{1}}, I_{B_{1}})$ and $(A_{B_{2}},
I_{B_{2}})$ are equivalent.

\sssec{Chevalley generators and Chevalley bases}\label{SsChev} We
often denote the set of generators corresponding to a normalized
matrix by $X_{1}^{\pm},\dots , X_{n}^{\pm}$ instead of
$e_{1}^{\pm},\dots , e_{n}^{\pm}$; and call them, together with the
elements $H_i:=[X_{i}^{+}, X_{i}^{-}]$ added for convenience for all
$i$, the {\it Chevalley generators}.\index{Chevalley generator}

For normalized Cartan matrices of simple finite dimensional Lie
algebras, there exists only one (up to signs) basis consisting of
vectors homogeneous with respect to the weight grading and
containing $X_i^\pm$ and $H_i$ in which all structure constants are
integer, cf. \cite{St}. Such a basis is called the\footnote{Observe
that, for a distinct normalization, there might exist another basis
with integer structure constants, and our normalization --- the most
natural, it seems, in super setting, cf. \cite{BGL} --- is
different, when specialized to $\fo(2n+1)$, from the Chevalley one.}
{\it Chevalley}\index{Basis! Chevalley} basis.

Observe that, having normalized the Cartan matrix of $\fo(2n+1)$ so
that $A_{nn}=1$, we get {\bf another} basis with integer structure
constants. We think that this basis also qualifies to be called {\it
Chevalley basis}; for Lie superalgebras such normalization is the
most reasonable.

\section{Cartan matrices and Dynkin
diagrams}

\ssec{Disclaimer} A usual way to represent simple Lie algebras over
$\Cee$ with integer Cartan matrices is via graphs called, in the
finite dimensional case, {\it Dynkin diagrams} (DD). The Cartan
matrices of certain interesting infinite dimensional simple Lie {\it
super}algebras $\fg$ over $\Cee$ can be non-symmetrizable or have
entries belonging to the ground field. Still, it is always possible
to assign an analog of the Dynkin diagram to each Lie (super)algebra
(with Cartan matrix, of course) provided the edges and nodes of the
graph (DD) are rigged with an extra information. Although these
analogs of the Dynkin graphs are not uniquely recovered from the
Cartan matrix (and the other way round), they give a graphic
presentation of the Cartan matrices and help to observe some hidden
symmetries.

Obviously, representation of Cartan matrices by means of DD is only
meaningful for rather sparse matrices. Cartan matrices of many of
the Lorentzian Lie (super)algebras have all their entries non-zero
and of large absolute value, making representation by means of DD
useless.

\ssec{Nodes} To every simple root there may correspond
\begin{equation}\label{cm1}
\begin{cases}
\text{a node}\; \mcirc\; &\text{if $p(\alpha_{i})= \ev$ and $
A_{ii}=2$},\\
\text{a node}\; \ast \;&\text{if $p(\alpha_{i}) =\ev$ and
$A_{ii}=\od$};\\
\text{a node}\; \mbullet \;&\text{if $p(\alpha_{i}) =\od$ and
$A_{ii}=1$};\\
\text{a node}\; \motimes \;& \text{if $p(\alpha_{i}) =\od$ and $
A_{ii}=0$},\\
\text{a node}\; \odot\; &\text{if $p(\alpha_{i})= \ev$ and $
A_{ii}=\ev$}.\\
\end{cases}
\end{equation}
In characteristic 0,  to construct Lie superalgebras of types Fin
and Aff with indecomposable Cartan matrices --- the puzzles of our
Lego problem to be described in what follows --- we do not need the
nodes $\ast$ and $\odot$ ({\it star} and {\it sun}, respectively).
The remaining nodes $\mcirc$, $\mbullet$ and $\motimes$ will be
referred to as {\it white}, {\it black} and {\it grey},
respectively.

The Lie algebra $\fsl(2)$ with Cartan matrix $(2)$, and the Lie
superalgebra $\fosp(1|2)$ with Cartan matrix $(1)$ are simple.

The Lie algebra with Cartan matrix $(\ev)$ and the Lie superalgebra
with Cartan matrix $(0)$ are solvable of $\dim 4$ and $\sdim 2|2$,
respectively. Their derived algebras --- {\it Heisenberg algebras}
--- are denoted $\fhei(2)\simeq\fhei(2|0)$ and $\fhei(0|2)\simeq\fsl(1|1)$,
respectively; their (super)dimensions are 3 and $1|2$,
respectively).

\ssec{Edges} If the Cartan matrix is integer, we can assign to it an
analog of the Dynkin diagram: Connect the $i$th node with the $j$th
one by $\max(|A_{ij}|, |A_{ji}|)$ edges rigged with an arrow $>$
pointing from the $i$th node to the $j$th if $|A_{ij}|>|A_{ji}|$ or
in the opposite direction if $|A_{ij}|<|A_{ji}|$. If the matrix
elements are not integer, we draw just one edge labeled with
$(|A_{ij}|, |A_{ji}|)$ as for $NS3_{82}-NS3_{86}$; for
$\fd(\alpha)$, we illustrate with an integer value of $\alpha$, but
with a label $\max(|A_{ij}|, |A_{ji}|)$ as in the generic case.

\ssec{Reflections} Let $R^+$ be a system of positive roots of Lie
superalgebra $\fg$, and let $B=\{\sigma_1,\dots,\sigma_n\}$ be the
corresponding system of simple roots with some corresponding pair
$(A=A_B,I=I_B)$. Then for any $k\in \{1, \dots, n\}$, the set
$(R^+\backslash\{\sigma_k\})\coprod\{-\sigma_k\}$ is a system of
positive roots. This operation is called {\it the reflection in
$\sigma_k$}; it changes the system of simple roots by the formulas
\begin{equation}
\label{oddrefl}
r_{\sigma_k}(\sigma_{j})= \begin{cases}{-\sigma_j}&\text{if~}k=j,\\
\sigma_j+B_{kj}\sigma_k&\text{if~}k\neq j,\end{cases}\end{equation}
where
\begin{equation}
\label{Boddrefl}B_{kj}=\begin{cases}
-\displaystyle\frac{2A_{kj}}{A_{kk}}& \text{~if~} A_{kk}\neq
0\text{~and~}
-\displaystyle\frac{2A_{kj}}{A_{kk}}\in \Zee_+, \\
1&\text{~if~}i_k=\od, A_{kk}=0,A_{kj}\neq 0,\\
0&\text{~if~}i_k=\od, A_{kk}=A_{kj}=0,\\
0&\text{~if~}i_k=\ev, A_{kk}=\ev,A_{kj}=0.\end{cases}
\end{equation}

If $A_{kk}\neq 0$ and $\displaystyle-\frac{2A_{kj}}{A_{kk}}\not\in
\Zee_+$, then the reflections are not defined (if the characteristic
of the ground field is equal to 0).

\sssbegin{Remark} The name \lq\lq reflection'' is used because in
the case of (semi)simple finite dimensional Lie algebras this
action, extended on the whole $R$ by linearity, is a map from $R$ to
$R$, and it does not depend on $R^+$, only on $\sigma_k$. This map
is usually denoted by $r_{\sigma_k}$ or just $r_{k}$. The map
$r_{\sigma_i}$ extended to the $\Ree$-span of $R$ is reflection in
the hyperplane orthogonal to $\sigma_i$ relative the bilinear form
dual to the Killing form.

The reflections in the even (odd) roots are referred to as {\it
even} ({\it odd}) {\it reflections}.\index{Reflection!
odd}\index{Reflection! even! non-isotropic} \index{Reflection! even!
isotropic} A simple root is called {\it isotropic}, if the
corresponding row of the Cartan matrix has zero on the diagonal, and
{\it non-isotropic} otherwise. The reflections that correspond to
isotropic or non-isotropic roots will be referred to accordingly.

If there are isotropic simple roots, the reflections $r_\alpha$ do
not, as a rule, generate a version of the {\it Weyl group} because
the reflections in a simple root is only defined for systems of
simple roots containing this root. In the general case, the action
of a given isotropic reflections (\ref{oddrefl}) can not, generally,
be extended to a linear map $R\tto R$. For Lie superalgebras over
$\Cee$, one can extend the action of reflections by linearity to the
root lattice but this extension preserves the root system only for
$\fsl(m|n)$
and $\fosp(2m+1|2n)$. 

If $\sigma_i$ is an odd isotropic root, then the corresponding
reflection sends one set of Chevalley generators into a new one:
\begin{equation}
\label{oddrefx} \tilde X_{i}^{\pm}=X_{i}^{\mp};\;\; \tilde
X_{j}^{\pm}=\begin{cases}[X_{i}^{\pm},
X_{j}^{\pm}]&\text{if $A_{ij}\neq 0$},\\
X_{j}^{\pm}&\text{otherwise}.\end{cases}
\end{equation}

\end{Remark}

\subsubsection{Serganova's lemma.}\label{Serle} Serganova \cite{Se, HS}
proved the following

\sssbegin{Lemma} For any Lie superalgebra of the form $\fg(A)$ with
an indecomposable $A$ and of polynomial growth, and any pair of
systems of simple roots $B_1$ and $B_2$, there exists a chain of
reflections connecting $B_1$ with some system of simple roots $B'_2$
equivalent (in the sense of definition \ref{EqSSR}) to either $B_2$
or $-B_2$. \end{Lemma}

\subsection{How to recover the Cartan matrix from the Dynkin diagram.}\label{cmfromDD}
Note that for $\fg=\fd(\alpha)$ or $\fg=\fd(\alpha)^{(1)}$, the
Cartan matrix is integer only if $\alpha\in\Zee$.

It so happens that for Lie superalgebras of certain types the
normalized Cartan matrix $(A_{ij})$ (and the set of parities
$I=\{i_1,\dots,i_n\}$ of the corresponding simple roots can be
practically uniquely, up to an equivalence, recovered from the
Dynkin graph. At least, for the finite dimensional and affine
Kac-Moody Lie algebras, such recovering is possible along the
following procedure (except for $\fg=\fd(\alpha)$,
$\fd(\alpha)^{(1)}$ and $\fsvect_\alpha(1|2)$ whose Dynkin diagrams
are drawn by special rules):

\begin{enumerate}

\item
If the $i$th and the $j$th vertices are connected by $k$ edges and
the arrow points at the $j$th vertex, then $|A_{ji}|=k$ and
$|A_{ij}|=1$.

\item
If the $i$th and the $j$th vertices are connected by $k$ edges
without arrows, then $|A_{ij}|=|A_{ji}|=k$.

\item
If the $j$th vertex is $\mbullet$, then $A_{jj}=1$ and $i_j=\od$.

\item[]
If the $j$th vertex is  $\motimes$, then $A_{jj}=0$ and $i_j=\od$.

\item[]
If the $j$th vertex is $\mcirc$, then $A_{jj}=2$ and $i_j=\ev$.

\item[]
If the $j$th vertex is $\ast$, then $A_{jj}=\od$ and $i_j=\ev$.

\item[]
If the $j$th vertex is $\odot$, then $A_{jj}=\ev$ and$i_j=\ev$.

\item If
$A_{ii}\ne0$, then $A_{ji}\le0$ for any $j\ne i$.

\item
If $A_{ii}=0$, then the $i$th row is recovered, up to a factor $-1$,
as follows:

\begin{enumerate}
\item[a)]
If $A_{ij}\ne0$ for precisely one $j$, the sign of~$A_{ij}$ can be
selected randomly;

\item[b)]
If $A_{ij_1},A_{ij_2}\ne0$ for precisely two indices $j_1$ è~$j_2$,
then $A_{ij_1}A_{ij_2}<0$;

\item[c)]
If $A_{ij_1},A_{ij_2},A_{ij_3}\ne0$ for precisely three distinct
indices $j_1$, $j_2$ and $j_3$, such that $A_{j_1j_2}\ne0$ and
$A_{j_1j_3}=A_{j_2j_3}=0$, then $A_{ij_1}A_{ij_2}>0$ and
$A_{ij_1}A_{ij_3}<0$.
\end{enumerate}
\end{enumerate}

These rules allow us to recover the Cartan matrix from the Dynkin
graph in the almost affine cases with the same uncertainty as in
affine and finite dimensional cases.

\subsection{Notation in Tables \ref{t17.0}--\ref{t17.5}.}\label{s17.1.12}

In Table \ref{t17.0} we show how we number the vertices in the
Dynkin graphs of finite dimensional and almost affine Lie
superalgebras. The symbol $\cdot$ stands for $\mcirc$ or $\motimes$.
The number $|v|$ is equal to the number of vertices of the Dynkin
graph, $ng$ is the number of \lq\lq grey" nodes $\motimes$ among the
nodes denoted by $\cdot$, and $png$ is the parity of $ng$.

Looking at the graphs we see that all of them, except for graphs for
$\fsl^{(1)}$ and $\fpsq^{(2)}$, can be situated on either one or
three horizontal levels. The maximal subgraph lying on the middle
level
\[
\xymatrix
{ \cdot\ar@{-}[r]&\dots&\cdot\ar@{-}[l] }
\]
(for $\fsl^{(1)}$ or $\fpsq^{(2)}$, the {\bf whole} graph) is said
to be the \textit{middle subgraph}.\index{Ñðåäíÿÿ ïîääèàãðàììà}

Let us split the middle subgraph into maximal non-intersecting
segments containing not more than one grey node $\motimes$ in such a
way that this node would be in the right end of each segment, for
example:

\begin{figure}[ht]\centering
\includegraphics[scale=.8]{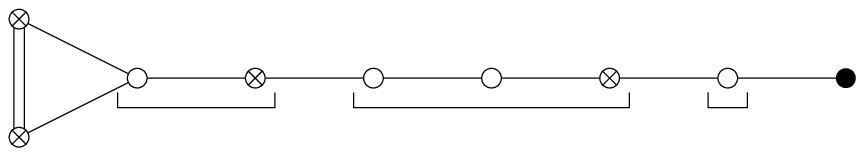}
\end{figure}
For the cyclic graph such a partition can be started anywhere, for
non-cyclic graphs, we begin from the left. Let us enumerate all
these segments, left to right or (for cyclic graphs)
counterclockwise. Let $ev$ (resp. $od$) be the total number of nodes
of the form $\motimes$ in all segments with even (odd) number of the
middle subdiagram.

All the Dynkin graphs of the types we consider can be uniquely
recovered from the data $v$, $ng$, $ev$, $od$.

A wavy line with arrows at its endpoints pointing at the grey
vertices  $\motimes$ of two diagrams in Tables signifies that these
diagrams are neighboring ones, connected by odd reflections in the
respective roots.

Unless otherwise stated, the numerical labels at vertices of the
graphs are coefficients of linear dependence of the rows of the
Cartan matrix.

\newpage\clearpage

\def\mmat #1,#2,#3,#4,{\text{\small$\displaystyle\left(\frac{#1}{#3}\frac{#2}{#4}\right)$}}
\newcommand{\smat}[1]{\text{\small\arraycolsep=2pt\def\arraystretch{1.1}%
           $\left(\begin{matrix}#1\end{matrix}\right)$}}
\def\tmat #1,#2,{\text{\small$\displaystyle\frac{#1}{#2}$}}

\rmnameii{Lim}{L}


\section[SRSes of types Fin and Aff]{Systems of simple roots of \lq\lq simple"
Lie superalgebras
of finite growth} \setcounter{table}{-1}

\begin{center}
\begin{surtab}

\tcaption{Skeletons and numbering of nodes}
\label{t17.0}\\[2mm]

\includegraphics[scale=.7]{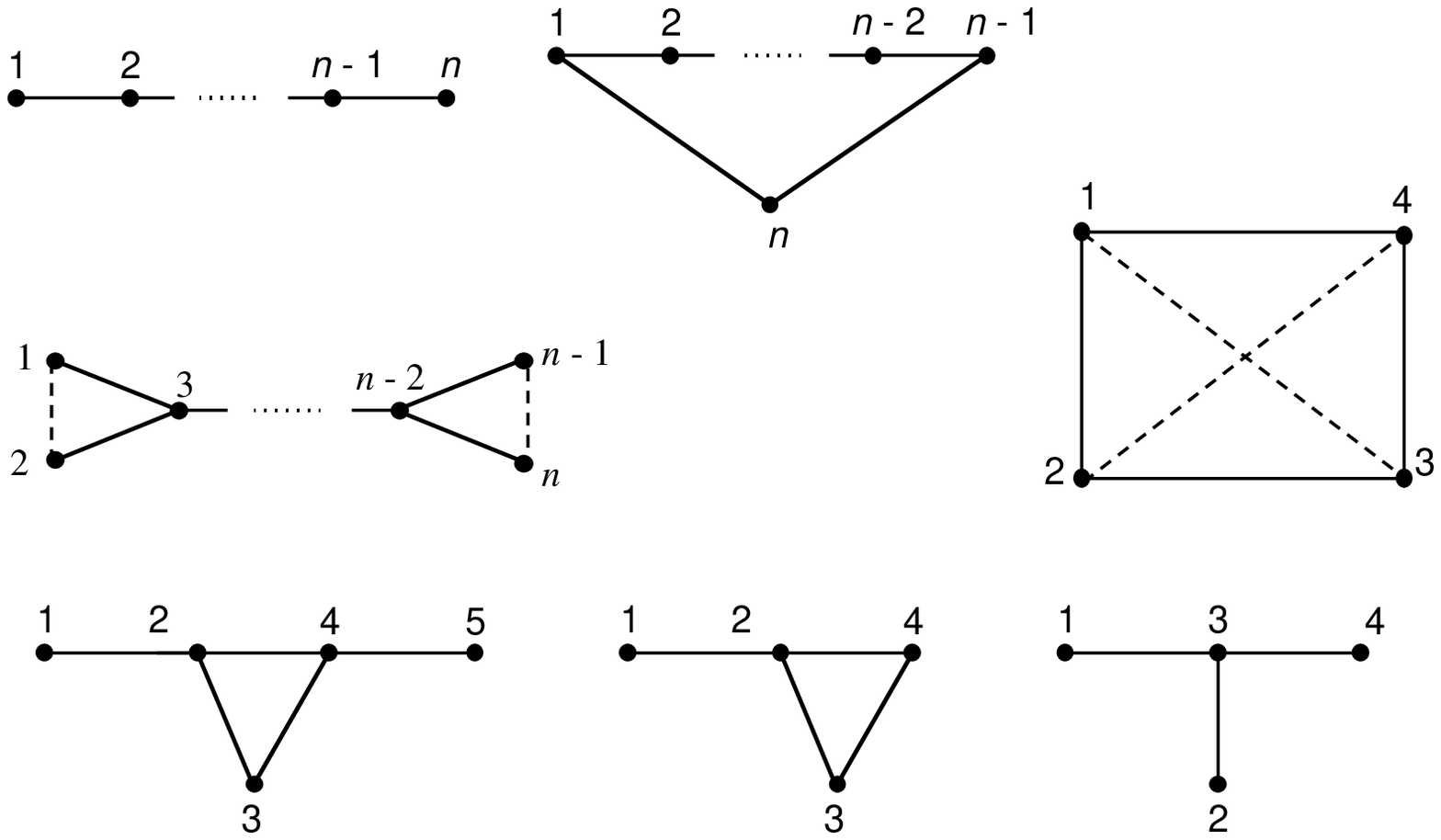}
\end{surtab}
\end{center}

\ssec{Non-symmetrizable Cartan matrices} There are the two series:

1) $\fsvect^L_\alpha(1|2)$: $
$  The odd reflections (in the second or third simple roots)
send $\alpha$ to $\alpha\pm 1$ (i.e., the above matrix turns into a
matrix of the same form but different value of parameter); so it
suffices to consider only the matrix depicted.

2) $\fpsq(n)^{(2)}$ (illustrated for $n=3$; $<-->$ stands for an odd
reflection):

$
$

\newpage

\begin{center}
\begin{surtab}

\setlength{\tabcolsep}{5pt}

\tcaption{Odd reflections}
\\[3mm] \extrarowheight=0pt \label{t17.2}

\begin{tabular}%
{>{\PBS\centering\hspace{0pt}}m{110mm}}
 \includegraphics[scale=1]{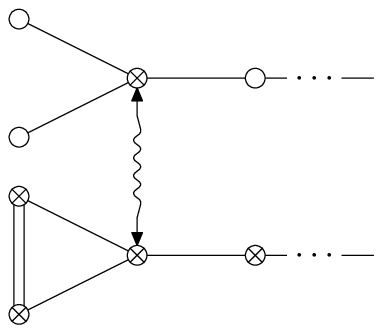}\quad
 \includegraphics[scale=1]{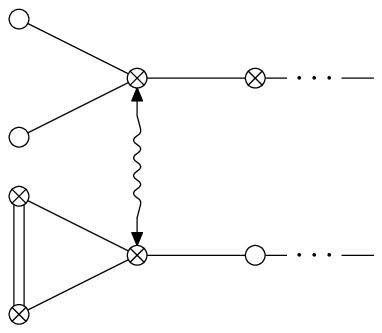}\\
 $png'=png+\od$,\; $od'=ev+1$,\; $ev'=od-1$\\[3mm]
 \includegraphics[scale=1]{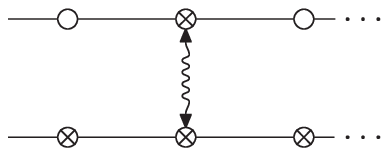}
 \includegraphics[scale=1]{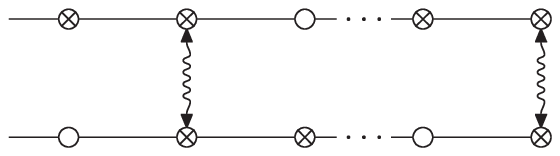}\\
 $\hskip 4.1em \underbrace{png'=png,\hskip 7.2em png'=png,\hskip 3.1em png'=png+\od}$\\
 $ev'=ev$,\; $od'=od$\\[3mm]
 \includegraphics[scale=1]{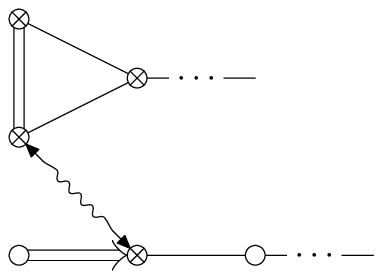}\quad
 \includegraphics[scale=1]{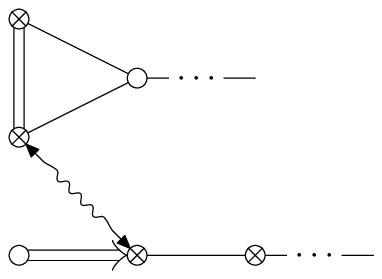}\\
 $png'=png$,\; $ev'=ev+1$,\; $od'=od$\\[3mm]
 \includegraphics[scale=1]{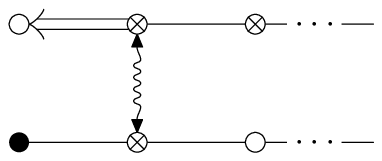}\quad
 \includegraphics[scale=1]{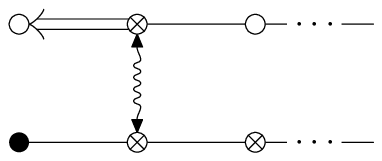}\\
 $png'=png+\od$,\; $od'=od$,\; $ev'=ev$\\[3mm]
\end{tabular}


\setlength{\tabcolsep}{0pt} \extrarowheight=0pt
\begin{tabular}
{>{\PBS\raggedright\hspace{0pt}}m{17mm}
 >{\PBS\raggedright\hspace{0pt}}m{40mm}
 >{\PBS\raggedright\hspace{0pt}}m{17mm}
 >{\PBS\raggedright\hspace{0pt}}m{40mm}}
 &&&\\[-5ex]
 $\fosp(4|2)^{(2)}$ && $\fsl(3|3)^{(4)}$&\\[-3ex]
 &\includegraphics{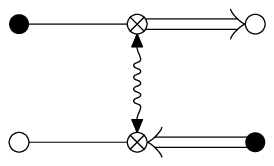}&& \includegraphics{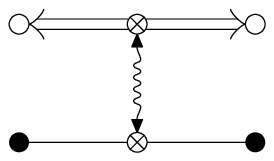}\\
 $\fsl(2|3)^{(2)}$ && $\fosp(3|2)^{(1)}$&\\[-3ex]
 &\includegraphics{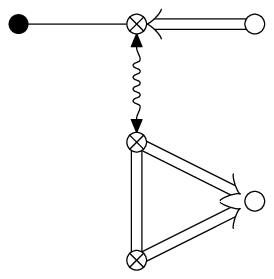}&& \includegraphics{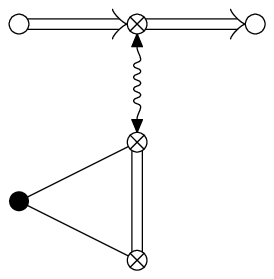}\\
\hline
\end{tabular}
\begin{tabular}%
{>{\PBS\raggedright\hspace{0pt}}p{57mm}
 >{\PBS\raggedright\hspace{0pt}}p{55.5mm}}
 $\fd(\alpha)$,\quad $\alpha=3$ \par
 \includegraphics[scale=.87]{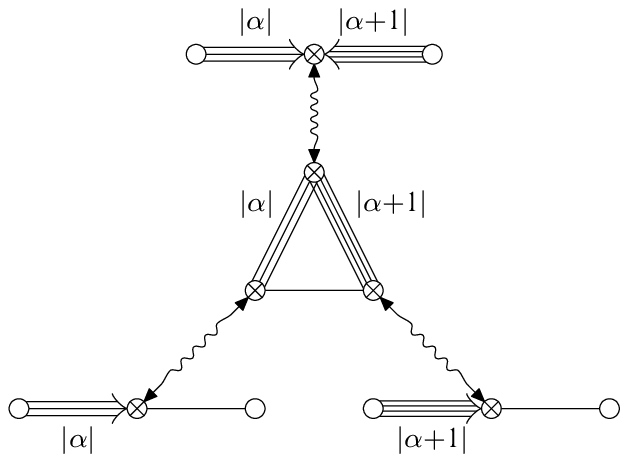}
 &
 $\fd(\alpha)^{(1)}$,\quad $\alpha=3$\par
 \includegraphics[scale=.87]{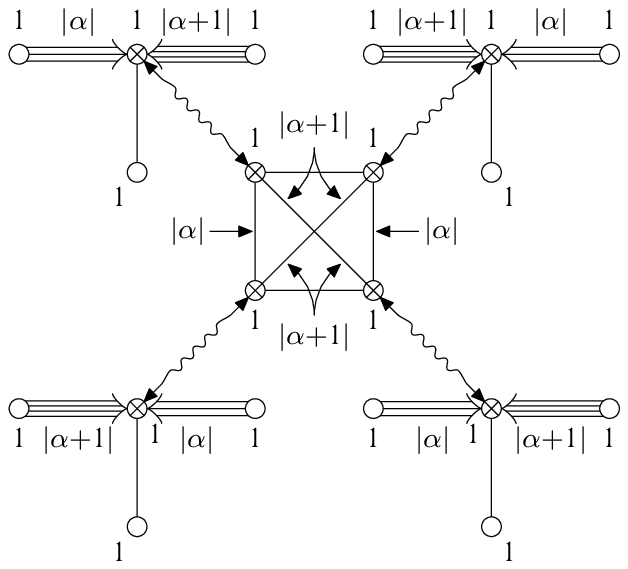}\\[-1mm]
\end{tabular}
\begin{tabular}%
{>{\PBS\raggedright\hspace{0pt}}p{38mm}
 >{\PBS\raggedright\hspace{0pt}}p{77mm}}
 $\fag(2)$\par
 \includegraphics[scale=1]{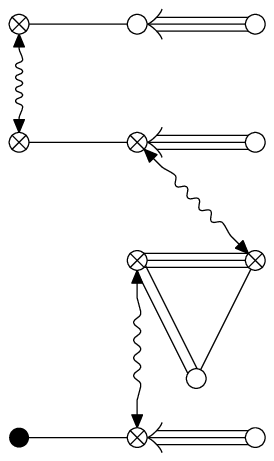}&
 $\fag(2)^{(1)}$\par
 \includegraphics[scale=1]{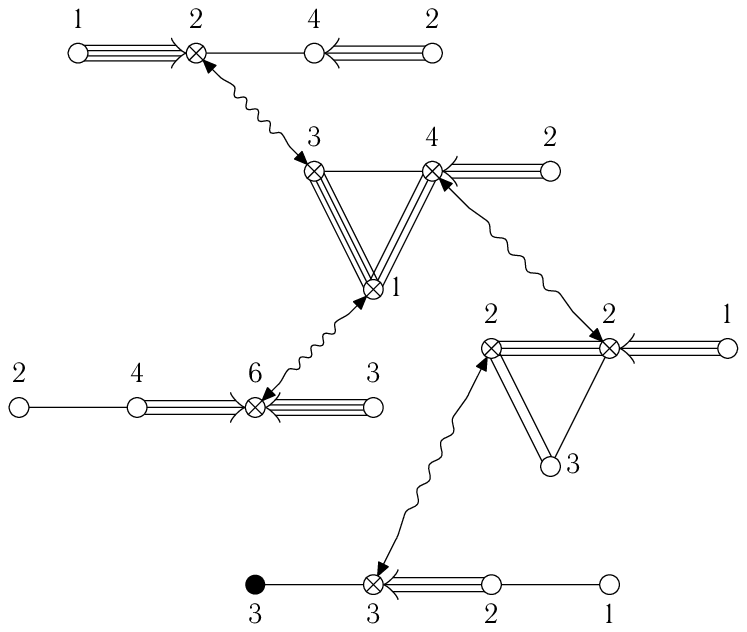}\\
\end{tabular}

\newpage
 \tabcolsep=4pt \extrarowheight=0pt
\begin{tabular}{ll}%
 $\fab(3)$&\\[-3ex]
 &\includegraphics[scale=1]{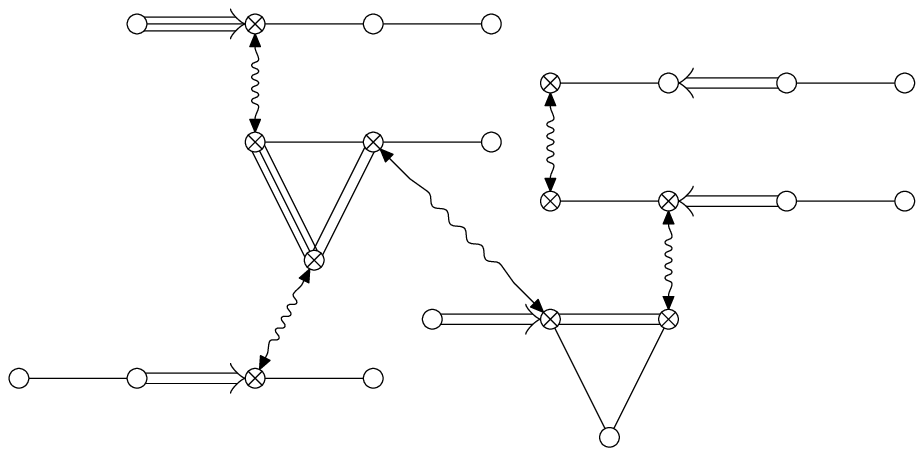}\\
\end{tabular}
\tabcolsep=0pt
\begin{tabular}{l}%
 $\fab(3)^{(1)}$\\
 \includegraphics[scale=.78]{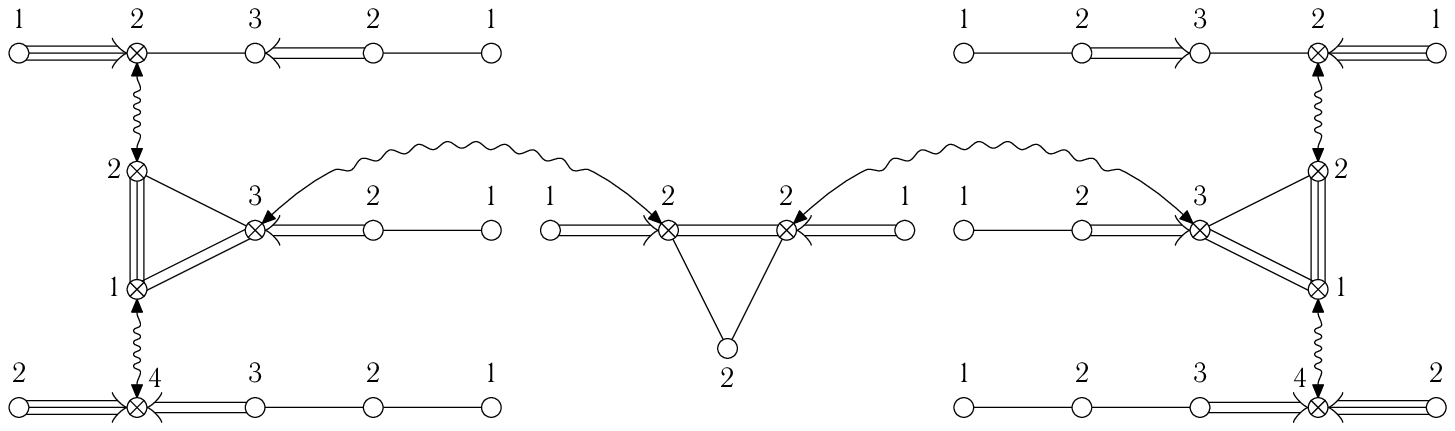}
\end{tabular}
\end{surtab}
\end{center}

\begin{landscape}
\thispagestyle{empty}
\begin{center}
\begin{surtab}
 \setlength{\tabcolsep}{0pt}
 \extrarowheight=0pt
 \tcaption{The simplest diagrams and their extensions}
 \\[0mm]
\label{t17.3}
\begin{tabular}%
{>{\PBS\raggedright\hspace{0pt}}m{22mm}
 >{\PBS\raggedright\hspace{0pt}}m{68mm}
 >{\PBS\raggedright\hspace{0pt}}m{70mm}
 >{\PBS\raggedleft\hspace{0pt}}m{26mm}}
 &&&$\fsl(n|m)^{(1)}$\\[-3ex]
 $\fsl(n|m)$,\par $m\ne n$\par\strut \par $\fsl(n|n)$\par \mbox{} &
 $\left.\arraycolsep=0pt
 \begin{array}{l}
 \includegraphics{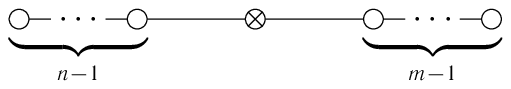}\\[1mm]
 \includegraphics{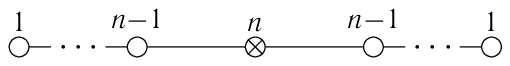}
 \end{array}\;
 \right\}$
 & \includegraphics{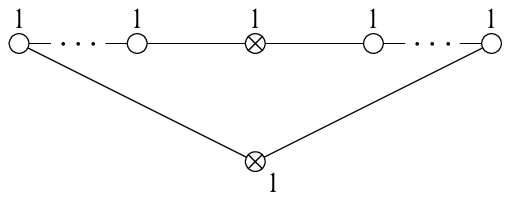}& \\[2ex]
 $\fosp(1|2n)$ &&&$\fosp(1|2n)^{(1)}$\\[-6ex]
 &\includegraphics{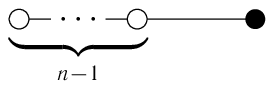} & \includegraphics{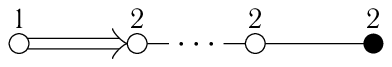} &\\[8ex]
 $\fosp(2|2n)$ &&& $\fosp(2|2n)^{(1)}$\\[-5ex]
 &\includegraphics{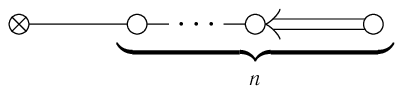} & \includegraphics{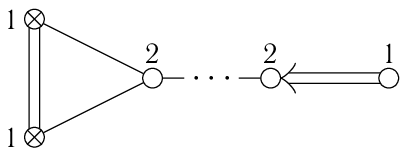} & \\[6ex]
 $\fosp(2m|2n)$ &&&\\[-4ex]
 &$\left.\arraycolsep=0pt
 \begin{array}{l}
 \includegraphics{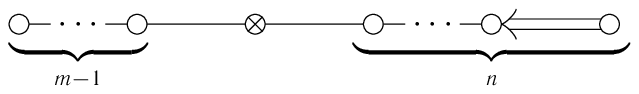}\\[1mm]
 \includegraphics{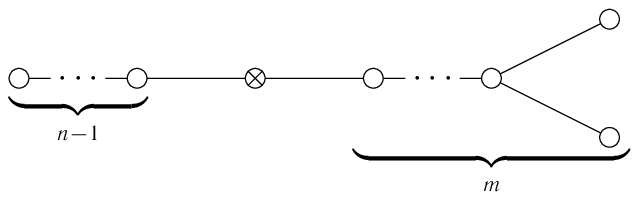}
 \end{array}\;
 \right\}\!\!$
 & \includegraphics{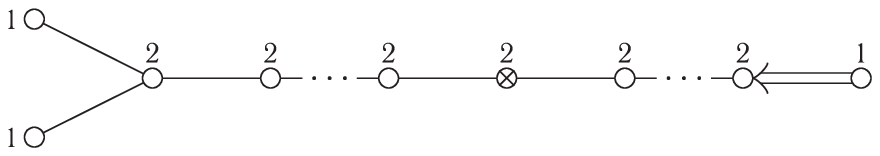} & \raisebox{8ex}[0pt]{$\fosp(2m|2n)^{(1)}$}\\
 $\fosp(2m+1|2n)$ &&\multicolumn{2}{r}{\raisebox{2.5ex}[0pt]{$\fosp(2m+1|2n)^{(1)}$}}\\[-2.5ex]
 &\includegraphics{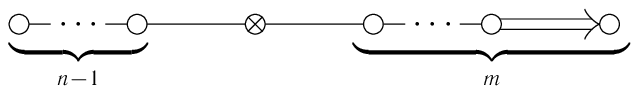} & \multicolumn{2}{l}{\includegraphics{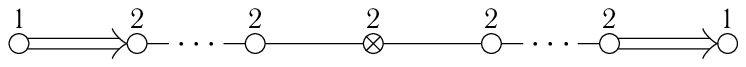}}  \\
\end{tabular}
\end{surtab}
\end{center}
\end{landscape}

\begin{center}
\begin{surtab}
 \tcaption{Systems of simple roots ($\rk\leq 4$)}
 \\[3mm]
\setlength{\tabcolsep}{0pt} \extrarowheight=0pt
\label{t17.4}
\begin{tabular}
{>{\PBS\raggedright\hspace{0pt}}m{22mm}
 >{\PBS\raggedright\hspace{0pt}}m{17mm}
 >{\PBS\raggedright\hspace{0pt}}m{18mm}|
 >{\PBS\raggedright\hspace{0pt}}m{22mm}
 >{\PBS\raggedright\hspace{0pt}}m{17mm}
 >{\PBS\raggedright\hspace{0pt}}m{16mm}}
 $\fg(A,I)$,\par $\sdim \fg(A,I)$ & \quad$A$ & Diagrams &
 \;\;$\fg(A,I)$,\par $\;\;\sdim \fg(A,I)$ & \qquad$A$ & Diagrams \\[3mm] \hline
 $\fsl(1|2)$,\par $4,4$ & $\smat{0&-1\\-1&2}$\par\mbox{}\par $\smat{0&1\\-1&0}$
 & \includegraphics{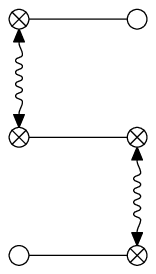} &
 \;\;$\fosp(3|2)$,\par $\;\;\sdim=(6,6)$ & $\smat{0&-1\\-2&2}$
 $\smat{
 0&1\\-1&1}$
 & \includegraphics{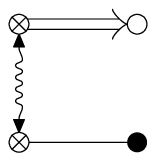}\\ 
 $\fosp(1|4)$,\par $\sdim=(10,4)$ & $\smat{2&-1\\-1&1}$
 & \includegraphics{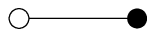} &
 \;\;$\fosp(1|2)^{(1)}$ & $\smat{2&-1\\-2&1}$
 & \includegraphics{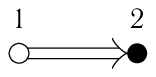}\\[2mm] 
 $\fosp(2|2)^{(2)}$ & $\smat{1&-1\\-1&1}$
 & \includegraphics{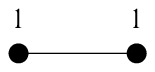} &
 \;\;$\fsl(1|3)^{(4)}$ & $\smat{2&-2\\-1&1}$
 & \includegraphics{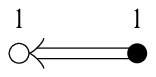}\\ 
\end{tabular}

\vskip 2mm

\setlength{\tabcolsep}{0pt} \extrarowheight=0pt
\begin{tabular}
{>{\PBS\raggedright\hspace{0pt}}m{23mm}
 >{\PBS\raggedright\hspace{0pt}}m{33mm}
 >{\PBS\raggedright\hspace{0pt}}m{23mm}
 >{\PBS\raggedright\hspace{0pt}}m{35mm}}
 &$\fsl(1|3)$, $\sdim=(9,6)$ && $\fpsl(2|2)$, $\sdim=(6,8)$ \\[-6mm]
 $\smat{0 &1 &0\\ -1&2&-1\\ 0&-1&2}$\par\mbox{}\par
 $\smat{0 &1 &0\\ -1&0&1\\ 0&-1&2}$
 &
 \vskip1em\includegraphics{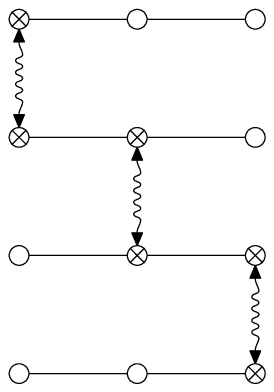}  & 
 $\smat{2 &-1 &0\\ -1&0&1\\ 0&-1&2}$\par\mbox{}\par
 $\smat{0 &1 &0\\ -1&0&1\\ 0&-1&0}$\par\mbox{}\par
 $\smat{0 &1 &0\\ -1&2&-1\\ 0&-1&0}$
 &\vskip1em\includegraphics{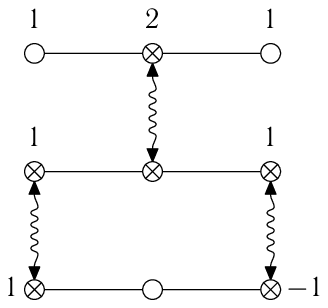}\vskip 3mm\par\mbox{}\\ 
 &&&\\[-4mm]
 $\smat{2 &-1 &0\\ -1&2&-1\\ 0&-1&1}$
 &  $\fosp(1|6)$,\par $\sdim=(21,6)$ \par\vskip 1em
 \includegraphics{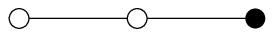}\vskip -1em &&\\
 &&&\\[-6mm]
 &$\fosp(3|4)$,\par $\sdim=(13,12)$ && $\fosp(5|2)$,\par $\sdim=(13,10)$ \\[-10mm]
 $\smat{2 &-1 &0\\ -1&0&1\\ 0&-2&2}$\par\mbox{}\par
 $\smat{0 &1 &0\\ -1&0&1\\ 0&-1&1}$\par\mbox{}\par
 $\smat{0 &1 &0\\ -1&2&-1\\ 0&-1&1}$
 & \vskip 3em\includegraphics{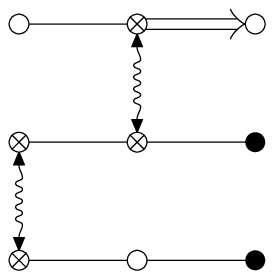}
 &$\smat{2 &-1 &0\\ -1&0&1\\ 0&-1&1}$\par\mbox{}\par
 $\smat{0 &1 &0\\ -1&0&1\\ 0&-2&2}$\par\mbox{}\par
 $\smat{0 &1 &0\\ -1&2&-1\\ 0&-2&2}$
 & \vskip 3em\includegraphics{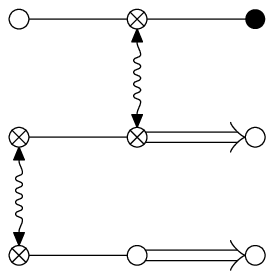}\\ 
 \end{tabular}

\begin{tabular}
{>{\PBS\centering\hspace{0pt}}m{40mm}
 >{\PBS\raggedright\hspace{0pt}}m{65mm}}
 &$\fd(\alpha)$,\; $\alpha=1$,\; $\sdim=(9,8)$ \\[-3em]
 $\smat{2 &-1 &0\\ \alpha&0&-1-\alpha\\ 0&-1&2}$\par
 $\smat{2 &-1 &0\\ -1&0&-\alpha\\ 0&-1&2}$\par
 $\smat{0 &1 &-1-\alpha\\ -1&0&-\alpha\\ -1-\alpha & \alpha&0}$\par
 $\smat{2 &-1 &0\\ -1-\alpha&0&1\\ 0&-1&2}$&
 \vskip3em\includegraphics[scale=.97]{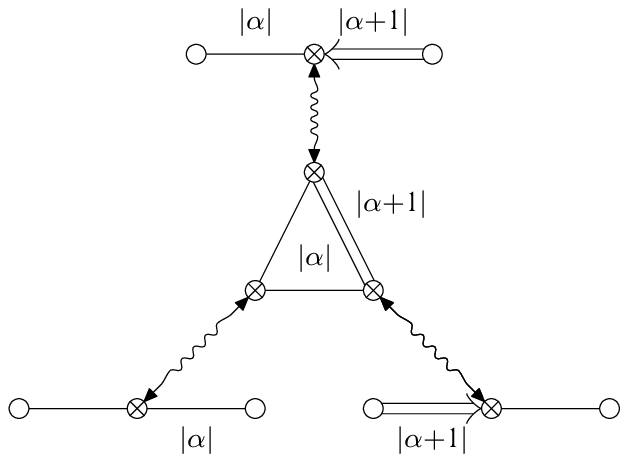}\\[0mm]
\end{tabular}
\begin{tabular}
{>{\PBS\raggedright\hspace{0pt}}m{22mm}
 >{\PBS\raggedright\hspace{0pt}}m{42mm}
 >{\PBS\raggedright\hspace{0pt}}m{22mm}
 >{\PBS\centering\hspace{0pt}}m{27mm}}
 &$\fag(2)$, $\sdim=(17,14)$&&\\[-.5cm]
 $\smat{0 &1 &0\\ -1&2&-3\\ 0&-1&2}$\par\mbox{}\par
 $\smat{0 &1 &0\\ -1&0&3\\ 0&-1&2}$\par\mbox{}\par
 $\smat{0 &-3 &1\\ -3&0&2\\ -1&-2&2}$\par\mbox{}\par
 $\smat{2 &-1 &0\\ -3&0&2\\ 0&-1&1}$
 & \mbox{}\par\vskip2mm\includegraphics{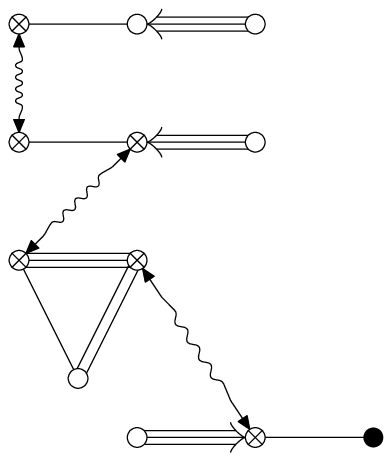}&
 \vskip-2.5cm $\smat{2 &-1 &-1\\ -1&0&1\\ -1&1&0}$
 & \vskip-2.8cm $\fsl(1|2)^{(2)}$\par
 \vskip1ex \includegraphics{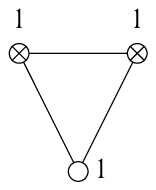}\vskip2cm\\[-5.5em] 
 && $\smat{2 &-1 &0\\ -2&2&-1\\ 0&-1&1}$
 & $\fsl(1|4)^{(1)}$\par \vskip1ex\includegraphics{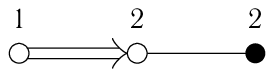}\\ 
\end{tabular}
\begin{tabular}
{>{\PBS\centering\hspace{0pt}}m{40mm}
 >{\PBS\raggedright\hspace{0pt}}m{65mm}}
 &\\[-1ex]
 &$\fosp(3|2)^{(1)}$ \\[-2em]
 $\smat{2 &-1 &0\\ 2&0&-1\\ 0&-2&2}$\par\mbox{}\par
 $\smat{0 &-2 &1\\ -2&0&1\\ -1&-1&1}$\par\mbox{}\par
 & \vskip2ex\includegraphics{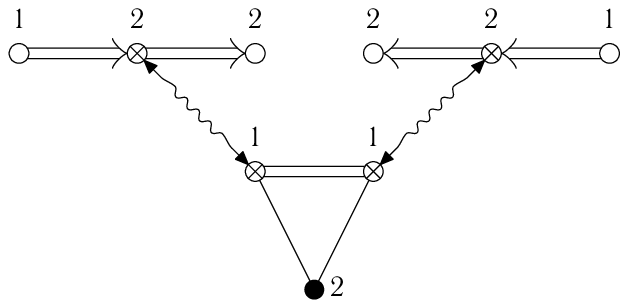}\\[3mm] 
\end{tabular}
\begin{tabular}
{>{\PBS\raggedright\hspace{0pt}}m{24mm}
 >{\PBS\raggedright\hspace{0pt}}m{35mm}
 >{\PBS\raggedright\hspace{0pt}}m{24mm}
 >{\PBS\raggedright\hspace{0pt}}m{30mm}}
 &$\fsl(1|4)^{(2)}$ &&  $\fosp(2|6)^{(2)}$ \\[-2em]
 $\smat{2 &-1 &0\\ -1&1&-1\\ 0&-1&2}$
 & \vskip2ex\includegraphics{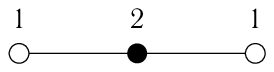}&
 $\smat{1 &-1 &0\\ -1&2&-1\\ 0&-1&1}$
 & \vskip2ex\includegraphics{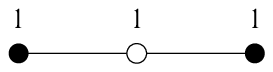}\\ 
\end{tabular}

\begin{tabular}
{>{\PBS\centering\hspace{0pt}}m{25mm}
 >{\PBS\centering\hspace{0pt}}m{65mm}}
 &$\fosp(4|2)^{(2)}$ \\[-2em]
 $\smat{1 &-1 &0\\ -1&0&1\\ 0&-2&2}$\par\mbox{}\par
 $\smat{2 &-2 &0\\ -2&0&2\\ 0&-1&1}$
 & \vskip1ex\includegraphics{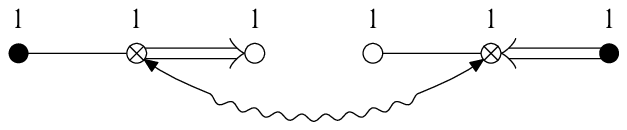}\\ 
 &\\[-3mm]
 &$\fsl(2|3)^{(2)}$ \\[-2em]
 $\smat{2 &-1 &0\\ -2&0&1\\ 0&-1&1}$\par\mbox{}\par
 $\smat{0 &2 &-1\\ -2&0&1\\ -2&-2&2}$\par\mbox{}\par
 $\smat{1 &-1 &0\\ -1&0&2\\ 0&-1&2}$
 & \vskip1ex\includegraphics{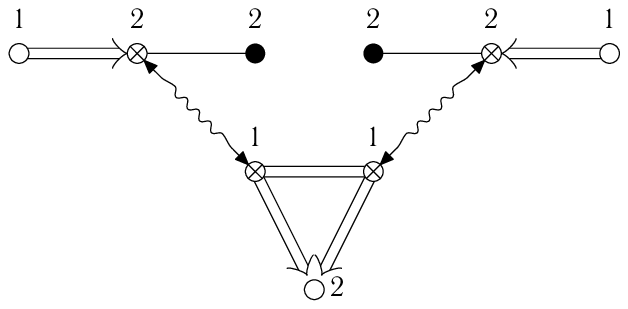}\\ 
\end{tabular}
\begin{tabular}
{>{\PBS\raggedright\hspace{0pt}}m{24mm}
 >{\PBS\raggedright\hspace{0pt}}m{35mm}
 >{\PBS\raggedright\hspace{0pt}}m{24mm}
 >{\PBS\raggedright\hspace{0pt}}m{30mm}}
 &&&\\[-3mm]
 $\smat{1 &-1 &0\\ -1&2&-1\\ 0&-2&2}$
 &  $\fsl(1|5)^{(4)}$ \par
 \vskip1ex\includegraphics{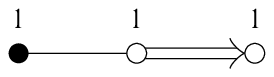}&&\\[-5.5em] 
 \vskip4.5em$\smat{2 &-2 &0\\ -1&0&1\\ 0&-2&2}$\par\mbox{}\par
 $\smat{1 &-1 &0\\ -1&0&1\\ 0&-1&1}$
 &  \vskip5.5em$\fsl(3|3)^{(4)}$\par
 \vskip1ex\includegraphics{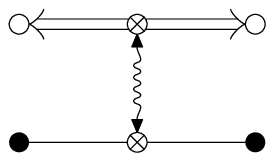}& 
 $\fosp(2|4)$,\; $\sdim=(11,8)$\par \vskip2em
 $\smat{0 &1 &0\\ -1&2&-2\\ 0&-1&2}$\par\mbox{}\par
 $\smat{0 &1 &0\\ -1&0&2\\ 0&-1&2}$\par\mbox{}\par
 $\smat{2&-1&-1\\1&0&-2\\-1&2&0}$
 & \vskip1pt\includegraphics{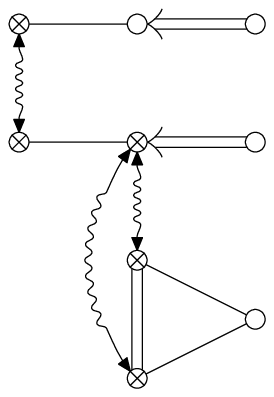}\\[1mm] 
\end{tabular}
\begin{tabular}
{>{\PBS\raggedright\hspace{0pt}}m{50mm}
 >{\PBS\raggedright\hspace{0pt}}m{40mm}}
 $\fsl(1|4)$,\; $\sdim=(16,8)$\par\vskip1ex
 \includegraphics{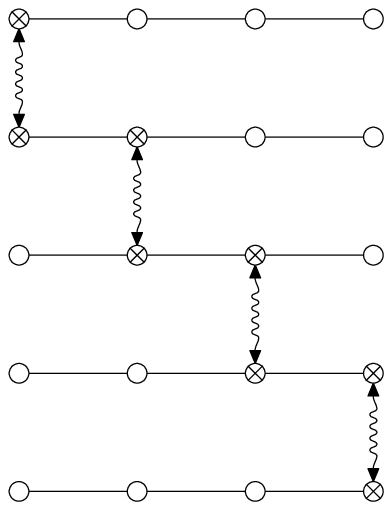}
 &
 $\fosp(1|8)$,\; $\sdim=(36,8)$\par\vskip1ex
 \includegraphics{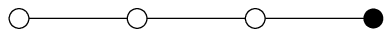}\par\vskip1ex
 $\fosp(3|6)$,\; $\sdim=(20,20)$\par\vskip1ex
 \includegraphics{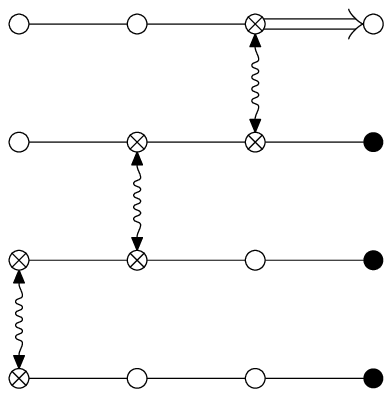}\\ 
\end{tabular}

\begin{tabular}
{>{\PBS\raggedright\hspace{0pt}}m{24mm}
 >{\PBS\raggedright\hspace{0pt}}m{88mm}}
 &\\[-9mm]
 $\fosp(5|4)$,\par $\sdim=(20,20)$ &\\[-3em]
 &\includegraphics{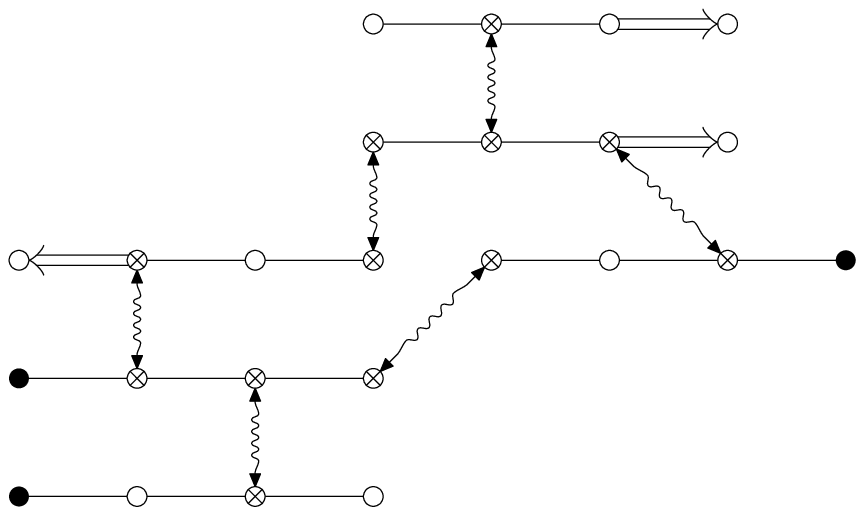}\\ 
 $\fsl(2|3)$,\par $\sdim=(12,12)$ &\\[-3em]
 &\includegraphics{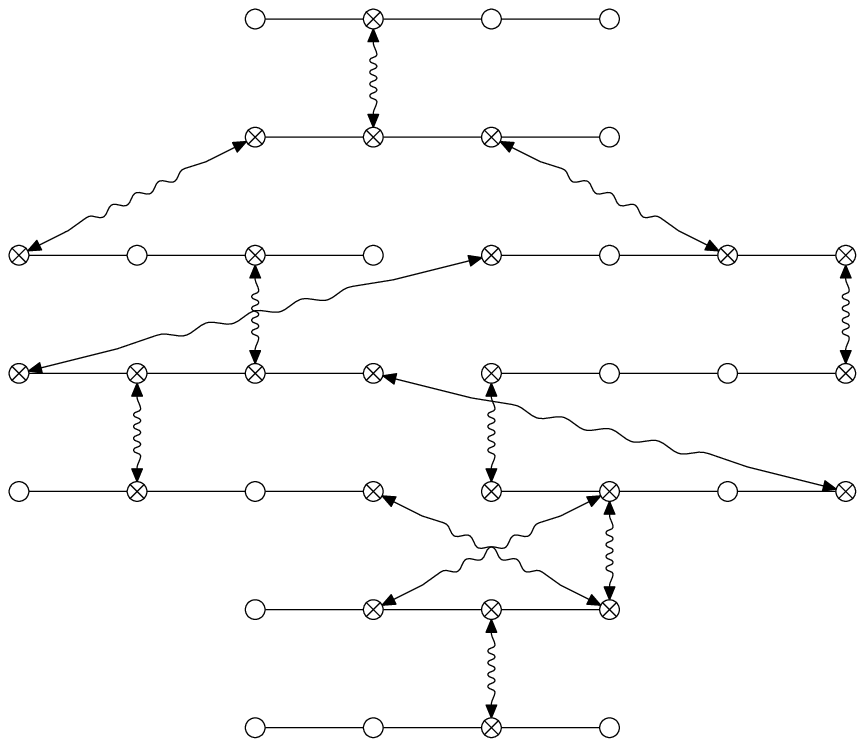}\\ 
 $\fosp(2|6)$,\par $\sdim=(24,12)$\\[-3em]
 &\includegraphics{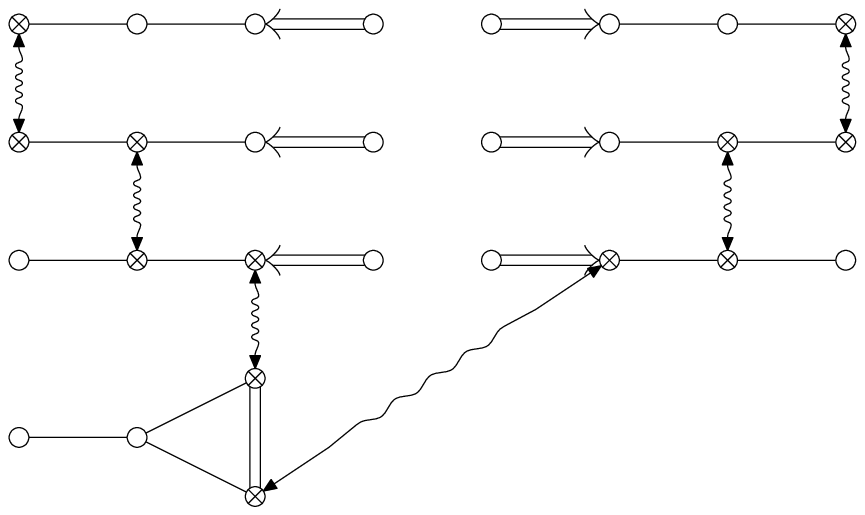}\\ 
\end{tabular}

\tabcolsep=3pt
\begin{tabular}
{>{\PBS\raggedright\hspace{0pt}}m{50mm}
 >{\PBS\raggedright\hspace{0pt}}m{50mm}}
 $\fosp(6|2)$,\raggedleft\par $\sdim=(18,12)$\raggedleft&\\[-3em]
 \vskip5em$\fosp(7|2)$,\par $\sdim=(24,14)$\par\vskip1ex
 \includegraphics{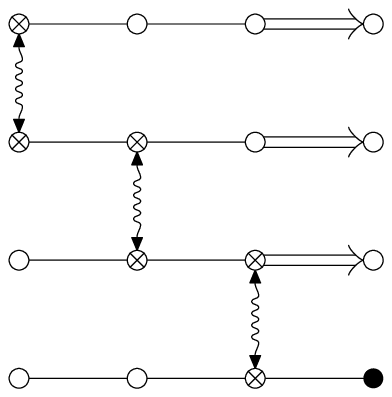}
 &
 \includegraphics{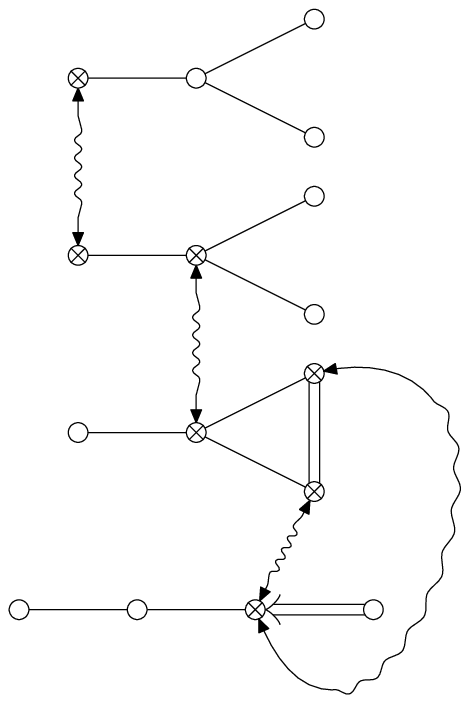}
 \\ 
\end{tabular}
\tabcolsep=0pt
\begin{tabular}
{>{\PBS\raggedright\hspace{0pt}}m{24mm}
 >{\PBS\raggedright\hspace{0pt}}m{85mm}}
 &\\[-4mm]
 $\fosp(4|4)$,\par $\sdim=(16,16)$ &\\[-3em]
 &\includegraphics{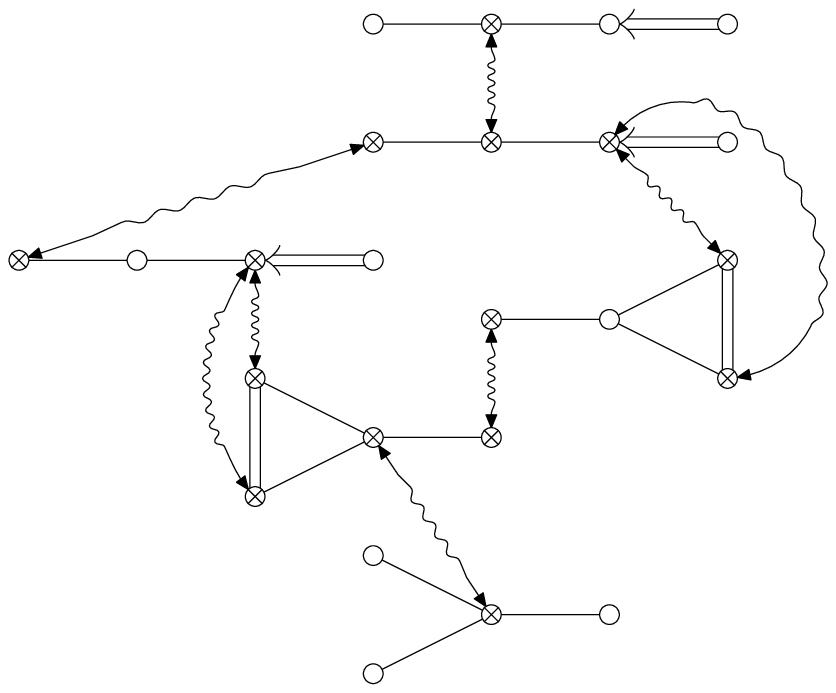}\\[2mm] 
\end{tabular}
\tabcolsep=0pt
\begin{tabular}
{>{\PBS\raggedright\hspace{0pt}}m{38mm}
 >{\PBS\raggedright\hspace{0pt}}m{65mm}}
 $\fsl(1|3)^{(1)}$&$\fsl(2|2)^{(1)}$\\
 \includegraphics{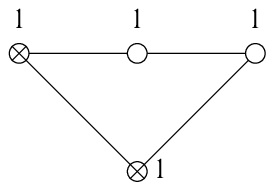}
 &\includegraphics{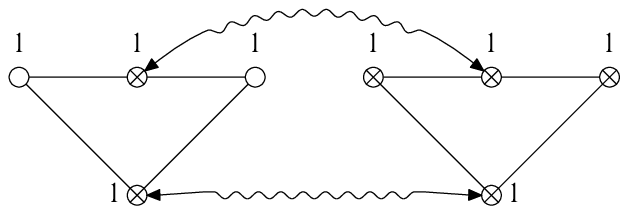}\\
\end{tabular}


\begin{tabular}
{>{\PBS\raggedright\hspace{0pt}}m{55mm}
 >{\PBS\raggedright\hspace{0pt}}m{56mm}}
 $\fosp(1|6)^{(1)}$\par\vskip1ex
 \includegraphics{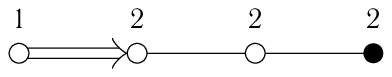}\par\vskip2ex
 $\fosp(3|4)^{(1)}$\par\vskip1ex
 \includegraphics{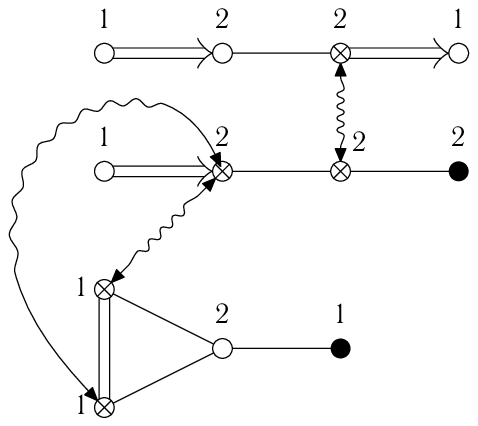}\par\vskip2ex
 $\fosp(2|4)^{(1)}$\par\vskip1ex
 \includegraphics{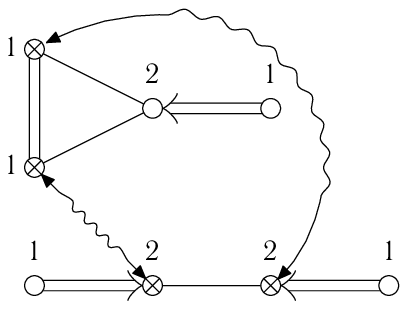}\par\vskip2ex
 $\fosp(2|4)^{(2)}$\par\vskip-.5em
 \hskip4.8em\includegraphics{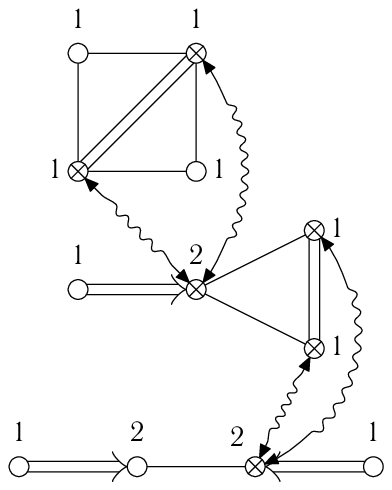}\par\vskip2ex
 $\fosp(2|6)^{(2)}$\par\vskip-.5em
 \hskip4.8em\includegraphics{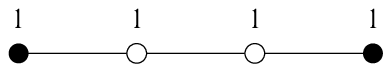}\par\vskip2ex
 &
 $\fosp(5|2)^{(1)}$\par\vskip.5em
 \hskip3.5em\includegraphics{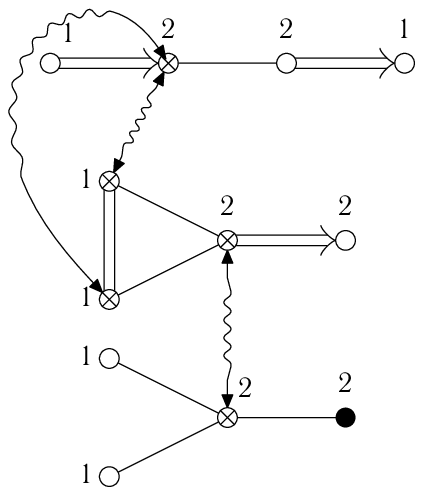}\par\vskip4ex
 $\fsl(1|6)^{(2)}$\par\vskip-1.5em
 \hskip4.5em\includegraphics{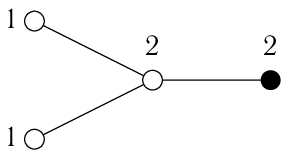}\par\vskip4ex
 $\fsl(2|5)^{(2)}$\par\vskip-1.5em
 \hskip4.5em\includegraphics{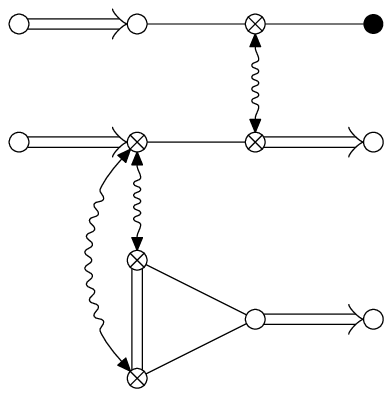}\par\vskip5ex
 $\fsl(3|4)^{(2)}$\par\vskip-2em
 \hskip4.5em\includegraphics{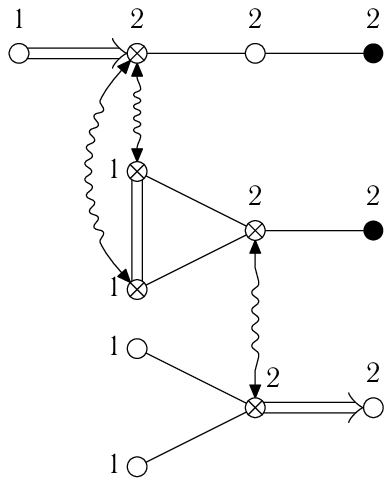}\\
\end{tabular}

\begin{tabular}
{>{\PBS\raggedright\hspace{0pt}}m{65mm}
 >{\PBS\raggedright\hspace{0pt}}m{45mm}}
 $\fosp(6|2)^{(2)}$\par\vskip-2em
 \hskip5em\includegraphics{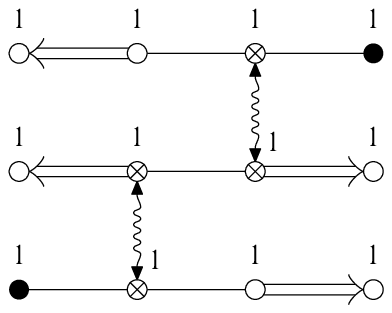}\par\vskip2ex
 $\fosp(4|4)^{(2)}$\par\vskip-2em
 \hskip5em\includegraphics{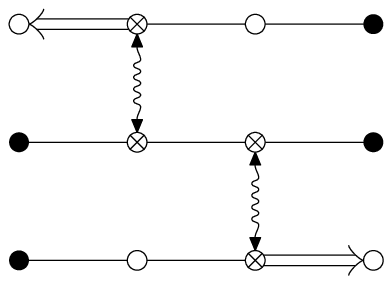}\par\vskip2ex
 &
 \vskip-5.5ex$\fsl(1|7)^{(4)}$\par\vskip1ex
 \includegraphics{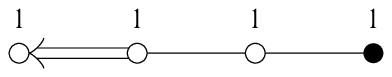}\par\vskip6ex
 $\fsl(3|5)^{(4)}$\par\vskip1ex
 \includegraphics{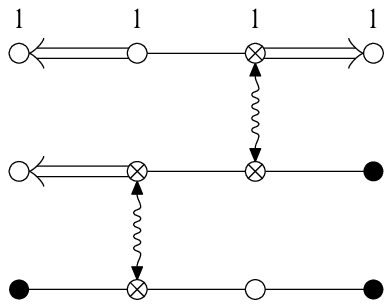}
\end{tabular}

\end{surtab}
\end{center}

\begin{landscape}
\begin{center}
\begin{surtab}
 \setlength{\tabcolsep}{0pt}
 \extrarowheight=0pt
 \tcaption{Systems of simple roots: General case}
 \\[3mm]
\label{t17.5}
\renewcommand{\thetable}{\arabic{chapter}.5.FD}
\begin{tabular}%
{>{\PBS\raggedright\hspace{0pt}}m{28mm}
 >{\PBS\raggedright\hspace{0pt}}m{31mm}
 >{\PBS\centering\hspace{0pt}}m{15mm}
 >{\PBS\centering\hspace{0pt}}m{15mm}
 >{\PBS\centering\hspace{0pt}}m{10mm}
 >{\PBS\centering\hspace{0pt}}m{8mm}
 >{\PBS\centering\hspace{0pt}}m{22mm}
 >{\PBS\centering\hspace{0pt}}m{20mm}}
 $\fg(A,I)$,\par $\sdim \fg(A,I)$ & Diagrams
 & $|v|$ & $ev\ge0$ & $od\ge0$ & $png$ & $ng\le \min$
 & Notation \`a la Cartan\\ \hline
 $\fsl(m|n)$,\; $m\ne n$\par $m^2+n^2-1,2mn$ &
 \includegraphics{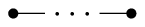} & $m+n-1$ & $m-1$ & $n$ & 1&
 $2m,2n$ & $A_{m,n-1}$\\
 $\fsl(n|n)$,\par $2n^2-2,2n^2$ &
 \includegraphics{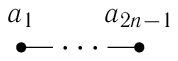} & $2n-1$ & $n-1$ & $n$ & &
 $2n-1$ & $A_{n,n-1}$\\
 $\fosp(1|2n)$,\par $2n^2+n,2n$ &
 \includegraphics{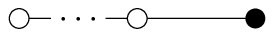} & $n$ & & & & & $B_{0,n}$\\
 $\fosp(2m|2n)$,\par $2n^2\,{+}\,2m^2\,{+}\,n\,{-}\,m,$ $4mn$ &
 $\left.\arraycolsep=0pt\begin{array}{l}
 \includegraphics{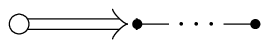}\\
 \includegraphics{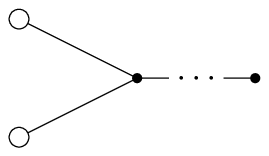}\\
 \includegraphics{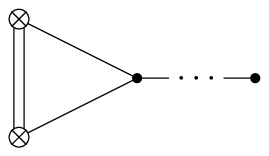}\\[-1ex]
 \end{array}\;\right\}$
 & $m+n$
 & $\arraycolsep=0pt\begin{array}{c}
   \tmat m,m-1,\\[2.2em] \tmat n,n-1,\\[3.2em] \tmat m-1,m-2,\\[2em]\end{array}$
 & $\arraycolsep=0pt\begin{array}{c}
   \tmat n-1,n,\\[2.2em] \tmat m-2,m-1,\\[3.2em] \tmat n-1,n,\\[2em]\end{array}$
 & $\arraycolsep=0pt\begin{array}{c}
   \tmat 0,1,\\[2.2em] \tmat 0,1,\\[3.2em] \tmat 0,1,\\[2em]\end{array}$
 & $\arraycolsep=0pt\begin{array}{c}
   2m,2n-1\\[2.2em] 2m-3,2n\\[3.2em] 2m-2,2n-1\\[2em]\end{array}$
 & $\arraycolsep=0pt\begin{array}{c}
   C_{m,n}\\[2.2em] D_{m,n}\\[3.2em] D^{\bcdot}_{m,n}\\[2em]\end{array}$\\
\hline
\end{tabular}
\end{surtab}
\end{center}
\end{landscape}

\begin{landscape}
\begin{center}
\begin{surtab}
 \setlength{\tabcolsep}{0pt}
 \extrarowheight=0pt
 \tcaption{Systems of simple roots: General case (cont.1)}
\\[3mm]
\begin{tabular}%
{>{\PBS\raggedright\hspace{0pt}}m{28mm}
 >{\PBS\raggedright\hspace{0pt}}m{47mm}
 >{\PBS\raggedright\hspace{0pt}}m{17mm}
 >{\PBS\centering\hspace{0pt}}m{14mm}
 >{\PBS\centering\hspace{0pt}}m{14mm}
 >{\PBS\centering\hspace{0pt}}m{10mm}
 >{\PBS\centering\hspace{0pt}}m{22mm}}
 $\fg(A,I)$,\par $\sdim \fg(A,I)$ & Diagrams & $|v|$\centering & $ev\ge0$ & $od\ge0$ & $png$ & $ng\le \min$\\ \hline
 &&&&&&\\[-3mm]
 $\fosp(2m+1|2n)^{(1)}$ &
 $\left.\arraycolsep=0pt\begin{array}{l}
 \includegraphics{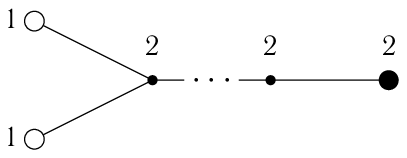}\\
 \includegraphics{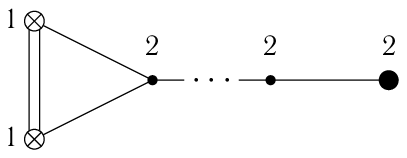}\\
 \includegraphics{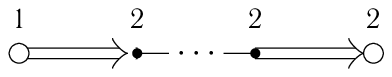}\\
 \includegraphics{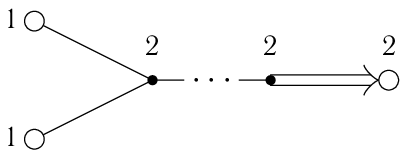}\\
 \includegraphics{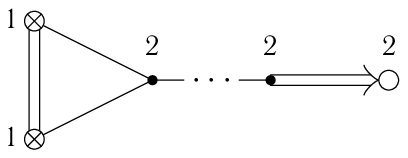}\\
 \includegraphics{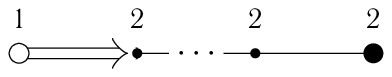}\\[-2mm]
 \end{array}\;\right\}$
 & $m+n+1$
 & $\arraycolsep=0pt\begin{array}{c}n-1\\[3.2em]m-1\\[2.8em]m-1\\[2.8em]n\\[3.1em]m-2\\[3em]m\\[-2em]\end{array}$
 & $\arraycolsep=0pt\begin{array}{c}m-1\\[3.2em]n-1\\[2.8em]n\\[2.8em]m-2\\[3.1em]n\\[3em]n-1\\[-2em]\end{array}$
 & $\arraycolsep=0pt\begin{array}{c}1\\[3.2em]0\\[2.8em]1\\[2.8em]0\\[3.1em]1\\[3em]0\\[-2em]\end{array}$
 & $\arraycolsep=0pt
 \begin{array}{l}2m-3,2n-1\\[3.2em]2m-1,2n-2\\[2.8em]2m-1,2n-1\\[2.8em]2m-4,2n\\[3.1em]2m-3,2n-1\\[3em]2m,2n-2\\[-2em]\end{array}$
 \\
\end{tabular}

\newpage
 \tcaption{Systems of simple roots: General case (cont.2)}
 \\[3mm]
\begin{tabular}%
{>{\PBS\raggedright\hspace{0pt}}m{19mm}
 >{\PBS\raggedright\hspace{0pt}}m{47mm}
 >{\PBS\centering\hspace{0pt}}m{15mm}
 >{\PBS\centering\hspace{0pt}}m{12mm}
 >{\PBS\centering\hspace{0pt}}m{12mm}
 >{\PBS\centering\hspace{0pt}}m{10mm}
 >{\PBS\centering\hspace{0pt}}m{22mm}}
 $\fosp(2m|2n)^{(1)}$ &
 $\left.\arraycolsep=0pt\begin{array}{l}
 \includegraphics{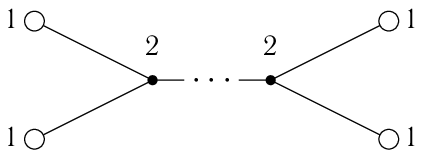}\\
 \includegraphics{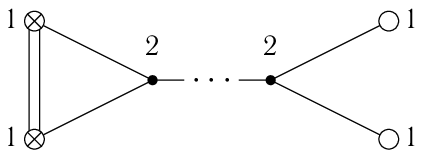}\\
 \includegraphics{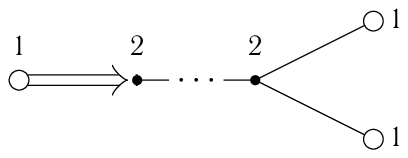}\\
 \includegraphics{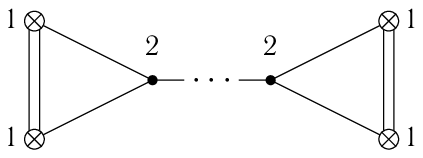}\\
 \includegraphics{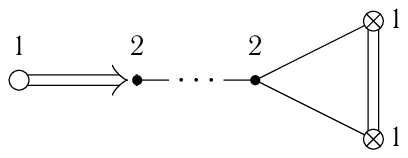}\\
 \includegraphics{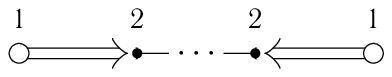}\\[-2mm]
 \end{array}\;\right\}$
 & $m+n+1$ & $\arraycolsep=0pt\begin{array}{c}n\\[3.8em]m-3\\[3.8em]m-2\\[3.8em]m-2\\[3.8em]m-1\\[3em]m\\[-2em]\end{array}$
 & $\arraycolsep=0pt\begin{array}{c}m-3\\[3.8em]n\\[3.8em]n\\[3.8em]n-1\\[3.8em]n-1\\[3em]n-1\\[-2em]\end{array}$
 & $\arraycolsep=0pt\begin{array}{c}0\\[3.8em]1\\[3.8em]1\\[3.8em]0\\[3.8em]0\\[3em]0\\[-2em]\end{array}$
 & $\arraycolsep=0pt
 \begin{array}{l}2m-6,2n\\[3.8em]2m-5,2n-1\\[3.8em]2m-3,2n-1\\[3.8em]2m-4,2n-2\\[3.8em]2m-2,2n-2\\[3em]2m,2n-2\\[-2em]\end{array}$\\
\end{tabular}

\newpage
 \tcaption{Systems of simple roots: General case (cont.3)}
 \\[1mm]
\begin{tabular}%
{>{\PBS\raggedright\hspace{0pt}}m{30mm}
 >{\PBS\raggedright\hspace{0pt}}m{47mm}
 >{\PBS\centering\hspace{0pt}}m{15mm}
 >{\PBS\centering\hspace{0pt}}m{12mm}
 >{\PBS\centering\hspace{0pt}}m{12mm}
 >{\PBS\centering\hspace{0pt}}m{10mm}
 >{\PBS\centering\hspace{0pt}}m{25mm}}
 $\fg(A,I)$,\par $\sdim \fg(A,I)$ & Diagrams & $|v|$ & $ev\ge0$ & $od\ge0$ & $png$ & $ng\le \min$\\ \hline
 &&&&&&\\[-5mm]
 $\fsl(m|n)^{(1)}$
 &\includegraphics{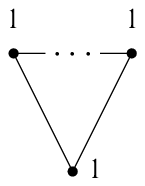} & $m+n$ & $\arraycolsep=0pt\begin{array}{l} m\\
 [-2mm] n\end{array}$
 & $\arraycolsep=0pt\begin{array}{l} n\\[-2mm] m\end{array}$ & 0 & $2m,2n$ \\
  \hline
$\fpsq(n)^{(2)}$
 &\includegraphics{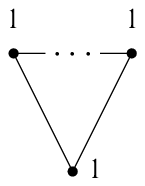} & $n$ & \multicolumn{2}{c}{$ev+od=n$} & 1&
 \\ \hline
 $\fosp(1|2n)^{(1)}$
 &\includegraphics{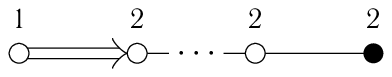} & $n+1$ & & & & \\
 $\fosp(2m|2n)^{(2)}$\par или\par $\fsl(2m+1|2n+1)^{(4)}_{-st}$ &
 $\arraycolsep=0pt\begin{array}{l}
 \includegraphics{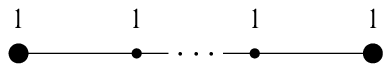}\\[1em]
 \includegraphics{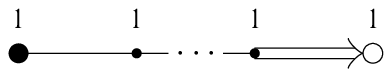}\\[1em]
 \includegraphics{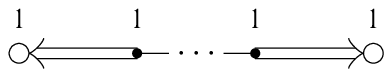}\\
 \end{array}$
 & $m+n$\par или\par $m+n+1$
 & $\arraycolsep=0pt\begin{array}{c}
    \tmat {m-1},{m},\\[1.99em]\tmat m-2,m,\\[1.99em]\tmat n,m-1,\\\end{array}$
 & $\arraycolsep=0pt\begin{array}{c}
    \tmat {n-1},{n-1},\\[1.99em]\tmat n,n-1,\\[1.99em]\tmat m-2,n,\\\end{array}$
 & $\arraycolsep=0pt\begin{array}{c}
    \tmat 0,1,\\[1.99em]\tmat 1,0,\\[1.99em]\tmat 0,1,\\\end{array}$
 & $\arraycolsep=0pt\begin{array}{l}
    \tmat {2m-2,2n-2},{2m-1,2n-1},\\[1.99em]\tmat {2m-3,2n-1},{2m,2n-2},\\[1.99em]
    \tmat {2m-4,2n},{2m-1,2n-1},\\ \end{array}$
 \\ \hline
\end{tabular}

\end{surtab}
\end{center}
\end{landscape}

\newpage
\begin{landscape}
\begin{center}
\begin{surtab}

\vspace{-3mm}
 \setlength{\tabcolsep}{0pt}
 \extrarowheight=0pt
 \tcaption{Systems of simple roots: General case (cont.4)}
 \\[3mm]
\begin{tabular}%
{>{\PBS\raggedright\hspace{0pt}}m{25mm}
 >{\PBS\raggedright\hspace{0pt}}m{47mm}
 >{\PBS\centering\hspace{0pt}}m{15mm}
 >{\PBS\centering\hspace{0pt}}m{12mm}
 >{\PBS\centering\hspace{0pt}}m{12mm}
 >{\PBS\centering\hspace{0pt}}m{10mm}
 >{\PBS\centering\hspace{0pt}}m{24mm}}
 $\fg(A,I)$,\par $\sdim \fg(A,I)$ & Diagrams & $|v|$ & $ev\ge0$ & $od\ge0$ & $png$ & $ng\le \min$\\ \hline
 &&&&&&\\[-3em]
 $\fosp(2m|2n)^{(2)}$ &
 $\left.\arraycolsep=0pt\begin{array}{l}
 \includegraphics{rrt5.11.eps}\\
 \includegraphics{rrt5.12.eps}\\
 \includegraphics{rrt5.13.eps}\\
 \includegraphics{rrt5.14.eps}\\
 \includegraphics{rrt5.15.eps}\\
 \includegraphics{rrt5.16.eps}\\[-2mm]
 \end{array}\;\right\}$
 & $m+n+1$
 & $\arraycolsep=0pt\begin{array}{c}
    n-2\\[3.5em]\tmat {m-1},{n-1},\\[3.8em]\tmat n,m,\\[3.8em]m-2\\[3.8em]\tmat n-1,m-1,\\[3em]n-1\\[-2em]\end{array}$
 & $\arraycolsep=0pt\begin{array}{c}
    m-1\\[3.5em]\tmat {n-2},{m-2},\\[3.8em]\tmat m-2,n-2,\\[3.8em]n-1\\[3.8em]\tmat m-1,n-1,\\[3em]m\\[-2em]\end{array}$
 & $\arraycolsep=0pt\begin{array}{c}
    1\\[3.8em]0\\[3.8em]0\\[3.8em]1\\[3.8em]1\\[3em]1\\[-2em]\end{array}$
 & $\arraycolsep=0pt\begin{array}{l}
    2m-3,2n-3\\[3.5em]\tmat {2m-2,2n-4},{2n-2,2m-4},\\[3.8em]\tmat {2m-4,2n},{2n-4,2m},\\[3.8em]
    2m-3,2n-3\\[3.8em]\tmat {2m-3,2n-1},{2n-3,2m-1},\\[3em]2m-1,2n-1\\[-2em]\end{array}$
\end{tabular}

\newpage
 \tcaption{Systems of simple roots: General case (cont.5)}
 \\[3mm]
\begin{tabular}%
{>{\PBS\raggedright\hspace{0pt}}m{25mm}
 >{\PBS\raggedright\hspace{0pt}}m{47mm}
 >{\PBS\centering\hspace{0pt}}m{15mm}
 >{\PBS\centering\hspace{0pt}}m{12mm}
 >{\PBS\centering\hspace{0pt}}m{12mm}
 >{\PBS\centering\hspace{0pt}}m{10mm}
 >{\PBS\centering\hspace{0pt}}m{22mm}}
$\fg(A,I)$,\par $\sdim \fg(A,I)$ & Diagrams & $|v|$ & $ev\ge0$ &
$od\ge0$ & $png$ & $ng\le \min$\\ \hline
 &&&&&&\\[-3em]
 $\fsl(2m+1|2n)^{(2)}_{-st}$ &
 $\left.\arraycolsep=0pt\begin{array}{l}
 \includegraphics{rrt5.17.eps}\\
 \includegraphics{rrt5.18.eps}\\
 \includegraphics{rrt5.19.eps}\\
 \includegraphics{rrt5.20.eps}\\
 \includegraphics{rrt5.21.eps}\\
 \includegraphics{rrt5.22.eps}\\[-2mm]
 \end{array}\;\right\}$
 & $m+n+1$
 & $\arraycolsep=0pt\begin{array}{c}
    m\\[3.2em]n-2\\[2.8em]\tmat n,n-1,\\[2.8em]m-1\\[3.1em]n-1\\[3em]n-1\\[-2em]\end{array}$
 & $\arraycolsep=0pt\begin{array}{c}
    n-2\\[3.2em]m\\[2.8em]\tmat m-1,m,\\[2.8em]n-1\\[3.1em]m-1\\[3em]m\\[-2em]\end{array}$
 & $\arraycolsep=0pt\begin{array}{c}
    0\\[3.2em]1\\[2.8em]\tmat 0,0,\\[2.8em]1\\[3.1em]0\\[3em]1\\[-2em]\end{array}$
 & $\arraycolsep=0pt\begin{array}{l}
  2m,2n-4\\[3.2em]2m-1,2n-3\\[2.8em]2m-2,2n\\[2.8em]2m-1,2n-3\\[3.1em]2m-2,2n-2\\[3em]2m,2n-2\\[-2em]\end{array}$\\
\end{tabular}

\end{surtab}
\end{center}
\end{landscape}

\newpage\clearpage

\section{Problems: Discussion}

Recently, Frappat and Sciarrino (\cite{FS}) offered a list of of
200-odd examples of what we call almost affine (they called them
\lq\lq hyperbolic") Lie superalgebras
\begin{equation}\label{intCM}
\text{with {\bf integer and symmetrizable} Cartan matrices }
\end{equation}
(correcting an earlier stake \cite{TDP} of classification). More
precisely, Frappat and Sciarrino attempted to classify \lq\lq
hyperbolic" (i.e., almost affine in our terms) {\bf Cartan matrices,
rather then the corresponding Lie superalgebras}.

In view of Remark \ref{tafrek}, (\ref{param}) and (\ref{non-sy}),
the restrictions \eqref{intCM} are manifestly unnatural. However,
even under these restrictions infinitely many examples (see Example
\ref{exa}) of \lq\lq hyperbolic" Cartan matrices if
$\text{size}(A)=3$ immediately spring to one's mind. Besides,
Frappat and Sciarrino \cite{FS} do not give a clear definition of
the Dynkin graphs (it is not needed, actually, since we are dealing
with Cartan matrices, but the result of \cite{FS} is given in terms
of Dynkin graphs, whatever they might be), neither is the strategy
of the search used in \cite{FS} clear (if there was any: Frappat and
Sciarrino expressed their hope for completeness of their result,
which means they can not offer means for verification to the
reader).

On top of all this, both \cite{TDP} and \cite{FS} deal with certain
(almost affine) properties of Cartan {\bf matrices} while speaking
about {\bf Lie superalgebras}. Unlike the non-super case, the
property of one Cartan matrix to be almost affine does not guarantee
same for the Cartan matrices obtained from this one by odd
reflection (however strange this might look for the inexperienced
reader).

All the above reasons combined, we decided to consider the following
Problems:

\ssbegin{Problem}\label{prob1} Classify almost affine Lie
superalgebras with indecomposable Cartan matrix with elements in the
ground field.
\end{Problem}

\sssbegin{Example}\label{exa} For any $a$ and $b$ (integer or not),
the Cartan matrix $\mat{0& 0& 1\\ 0& 0& 1\\ a& b& 0}$ is almost
affine, but
being reflected in the third simple root is turns into $\mat{2&\frac{a+b}{a}& -1\\
\frac{a+b}{b}&  2& -1\\ a & b &0}$ which is never almost affine,
whatever $a$ and $b$ (except for $a=-b$ in which case both matrices
are not almost affine because they correspond to $\fgl(2|2)$, or
--- by the usual abuse of language attributing Cartan matrices to simple Lie
algebras without Cartan matrices --- $\fpsl(2|2)$). This example
shows that the definition of almost affine Lie {\bf super}algebras
must differ from that in the non-super case: The \lq\lq almost
affine" property characterizes the given Cartan matrix rather than
the Lie {\bf super}algebra defined by means of this matrix.
\end{Example}

Problem \ref{prob1} is a typical Lego problem (to construct Lie
superalgebras  with indecomposable Cartan matrices from Lie
superalgebras of types Fin and Aff). It seems to be the easiest
among the problems listed in this section; at least, we have solved
it.


\ssec{Lorentzian Lie (super)algebras} Remarkable results of
Borcherds, Gritsenko and Nikulin require elucidating the following
relation.

\sssbegin{Problem}\label{prob2} Lorentzian Lie (super)algebras are
close to almost affine Lie (super)algebras; for rank small, they
sometimes coincide. Describe a precise relation between these types
of Lie superalgebras.\end{Problem}

Nikulin got some finiteness results concerning Lorentzian Lie
(super)algebras. It is not clear, however, if the problem is tame:
infinite series might be easier to describe than count grains of
sand on the beach, cf. \cite{N2} and \cite{VSh}.

\sssbegin{Problem}\label{probN} Classify the Lorentzian Lie
(super)algebras if the problem is tame.
\end{Problem}

One more incentive for our study: In our experience, various objects
(notions) of Linear Algebra and Calculus have \so{several}
counterparts in \lq\lq super" setting. So we wonder why each
Lorentzian Lie algebra $\fg_\ev$, considered so far, has at the
moment at most one almost affine Lie superalgebra $\fg$ with
$\fg_\ev$ as its even part.

\sssbegin{Problem}\label{probL} Which properties lack the \so{other}
almost affine Lie superalgebras having the same Lorentzian Lie
algebra as its even part?
\end{Problem}

\ssec{Deformations} Looking at simple finite dimensional and growth
1 Lie {\bf super}algebras some of which admit deformations
\so{preserving the same root system} it is natural to investigate
the following problem even in the non-super setting:

\sssbegin{Problem}\label{prob3} Are there deformations of Lorentzian
Lie (super)algebras? Are there deformations of almost affine Lie
(super)algebras?
\end{Problem}

Looking at \eqref{hysize2} it is clear that if $\rk\fg(A)=2$, there
are deformations even in the class of Lie algebras with Cartan
matrix. These deformations are, it seems, additional to those in the
class of Krichever-Novikov algebras. Amazingly, no complete
classification of deformations of affine Lie algebras are obtained
yet.

\ssec{Real forms} Having considered  almost affine Lie
(super)algebras over $\Cee$, we arrive at the following

\sssbegin{Problem}\label{prob4} Describe real forms of almost affine
Lie (super)algebras over $\Cee$.
\end{Problem}

\ssec{Presentations} A natural problem tackled in \cite{FS} is that
of explicit presentations:

\sssbegin{Problem}\label{prob5} Describe defining relations of
almost affine Lie (super)algebras (in terms of Chevalley
generators).
\end{Problem}

The answer is clear (all relations are of Serre type, even for
non-symmetrizable Cartan matrices), except for the cases
$NS3_{82}-NS3_{86}$. In these cases, for the answer, see
\cite{BGLr}.

\ssec{On equivalence classes of systems of simple roots} The point
is, for some infinite dimensional Lie (super)algebras with Cartan
matrix, there are several Cartan matrices related with each other by
means distinct from reflections. For example, the Lie algebras
$\fsl(\infty)$ whose Cartan matrix is infinite in all directions and
$\fsl_+(\infty)$ whose Cartan matrix goes only to the right and
downwards (or to the left and upwards) are isomorphic.

On the other hand, it is very interesting to find out
\begin{equation}\label{glinfty}
\begin{minipage}[l]{14cm}
{\sl What are the \lq\lq right" analogs of the Weyl group
(collection of reflections) that for a suitably enlarged Lie
superalgebra, there is only one equivalence class of systems of
simple roots?}
\end{minipage}\end{equation} For Egorov's approach to super versions
of $\fgl(\infty)$ and his solution, see \cite{Eg}.

\sssbegin{Problem}\label{prob6} Is an analog of Serganova's Lemma
\ref{Serle} true for almost affine Lie (super)algeb\-ras? In other
words: is there just one chain connecting two systems of simple
roots (perhaps, up to a change of sign)? Do we have to apply means
similar to those Egorov used? Same questions for Lorentzian Lie
(super)algebras.
\end{Problem}

\ssec{A purely super problem} Our result (Theorem \ref{theor}), as
well as classifications of the finite dimensional and (twisted) loop
algebras with indecomposable Cartan matrices, lead us to the
following problem indigenous to superalgebras.

\sssbegin{Problem} List all Lie superalgebras with indecomposable
Cartan matrix with at least one $0$ on the diagonal having
describable (either finitely many, or \lq\lq of the same type", as
for $\fsvect(1|2;a)$) set or sets of inequivalent Cartan matrices
connected by chains of odd reflections. As stated, this problem is
wild (or rather with a dull answer); single out its reasonable part,
if any.
\end{Problem}

\section{The algorithms}

\ssec{Skeletons and puzzle pieces}

1) In our supply of Lego puzzles, replace each CM by the one with
only $-1$ for every nonzero off-diagonal element and disregard the
difference between the nodes. Such CMs can be decoded by DDs that
will be referred to as {\it skeletons}. The extended skeletons
obtained by adding a node and drawing an edge from it to different
nodes of the given skeleton may correspond to different CMs, so the
nodes must be numbered, the new (added) node being the 0-th. Let the
number of nodes in a skeleton be its {\it size}.

2) Extending a given skeleton by induction on the size:

2.1) The new node can be joined with the 1-node skeleton.

2.2) The new node can be joined with any one of the two nodes of the
2-node skeleton, or with both the nodes.

2.3) The new node can be joined with any amount of nodes of any
3-node skeletons.

2.4) And so on. Beginning with 4-node  skeletons we encounter
impossible configurations.

2.end) Restriction on the number of nodes from above: impossibility
to extend any given skeleton of certain size and larger.

3) Recall that the nodes can be of different types (\lq\lq colors")
and list all possible colorations of the skeletons obtained at step
2).

4) For the colored skeletons obtained at step 3), consider possible
off-diagonal elements, or, equivalently, possibilities for multiple
edges and their directions.

5) Now recall that we are working with CMs, not DDs. By means of odd
reflections collect all CMs obtained at step 4) belonging to the
same Lie superalgebra, i.e., unite the CMs in orbits under \lq\lq
chains of odd reflections" that replace Weyl groups in super
setting.

This is easier to say than do, so we used computer to aid us; for
the program (without interface and practically without documentation
since we do not see in which situation to reuse this program), see
\cite{ChD}. We first tested the program by reproducing the result of
Li Wang Lai and manually checking the cases at variance with
\cite{FS}.

\ssec{Reflections} Let $A$ be a Cartan matrix of size $n$ and
$I=(p_1,...,p_n)$ the vector of parities. If $p_k=\od$ and
$A_{kk}=0$, then the reflection (\ref{oddrefl}) in the $k$th simple
odd root sends $A$ to
$$A'_{ij}=A_{ij} + b_iA_{kj} + c_jA_{ik},
$$where
$b_i=\begin{cases} -2&\text{if $i=k$;}\\
    0&\text{if $i\neq k$ and  $A_{ik}=0$;}\\
    \frac{A_{ik}}{A_{ki}}&\text{if $i\neq k$ and  $A_{ik}\neq 0$};\end{cases}$
and where
$c_j=\begin{cases} -2&\text{if $j=k$};\\
    0&\text{if $j\neq k$ and  $A_{jk}=0$};\\
    1&\text{if $j\neq k$ and  $A_{jk}\neq 0$}.\end{cases}$

After the reflection the matrix $A'$ might require to be normalized;
the new parities are $$p'_i=p_i+c_i \pmod 2.$$

This can be expressed in terms of matrices: $$A'=(E+B)A(E+C),$$
where all columns of the matrix $B$, except the $k$th one, are zero,
whereas the $i$th coordinate of the $k$th column-vector is $b_i$,
the $i$th coordinate of the $k$th row-vector of $C$ is $c_i$, the
other rows of $C$ being zero; $E$ is the unit matrix.

Serganova's lemma ensures us that, passing from one system of simple
roots (SSR) to a neighboring one, one can reach any other SSR (or
its opposite). To list all systems of simple roots is a standard
problem of graph chasing. Since the number of possible systems is
finite (for finite dimensional or affine Lie superalgebras of finite
rank) and not large, it does not take long. The algorithm is
standard: starting with a node, we first list all its neighbors,
then neighbors of its neighbors (neighbors of rank 2) except the
nodes already counted, and so on.

\section{Solution of Lego
problem} In this section, all parametric values are $\neq 0$.

\ssec{Size 2} We identify Cartan matrices obtained by permutation
of their rows and same columns; for $a,b\in \Ree$, we assume that
$a\leq b$.

Clearly, almost any (indecomposable) size 2 Cartan matrix yields an
almost affine Lie (super)algebra, namely there are the following
types:
\begin{equation}\label{size2}
\renewcommand{\arraystretch}{1.4}
\begin{array}{l}
1)\;\smat{2&a\\b&2},\quad
2)\;\smat{2&a\\b&1},\quad
3)\;\smat{2&a\\-1&0},\quad
4)\;\smat{1&a\\b&1},\quad
5)\;\smat{1&a\\-1&0},\quad
6)\;\smat{0&-1\\ -1&0}.\end{array}
\end{equation}
except case 6) which is of type Fin and the several other cases (easy to see)
which yield rank 2 Lie (super)algebras of types Fin or Aff.

In what follows, the Cartan matrices of almost affine Lie
superalgebras are listed. To the right of the Cartan matrix are
listed the finite dimensional or affine Lie algebras obtained after
deleting the respective row and column. The Lie algebras are listed,
for brevity, in Cartan's notations ($A_n=\fsl(n+1)$, etc.). If
deleting a vertex we get a disjoint union of graphs, one (or more)
of which is just a vertex, we do not indicate (to save space) the
obvious summands --- Lie (super)algebras $A_1$ and $\fosp(1|2)$.

\ssec{Size 3}{}~{}


\begin{tabular}
{>{\PBS\raggedright\hspace{0pt}}m{12mm}
 >{\PBS\raggedright\hspace{0pt}}m{27mm}
 >{\PBS\raggedright\hspace{0pt}}m{35mm}
 >{\PBS\raggedright\hspace{0pt}}m{30mm}}
 &&&\\
 $NS3_1)$&
 $\smat{2&-1&-1\\-1&1&-1\\
    -2&-1&2}$,&
 $\smat{1:\fosp(1|4)\\
    2:B_2\\
    3:\fosp(1|4)};
 $&
 \includegraphics{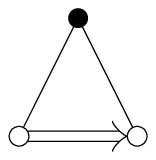}
\end{tabular}

\begin{tabular}
{>{\PBS\raggedright\hspace{0pt}}m{12mm}
 >{\PBS\raggedright\hspace{0pt}}m{27mm}
 >{\PBS\raggedright\hspace{0pt}}m{35mm}
 >{\PBS\raggedright\hspace{0pt}}m{30mm}}
 &&&\\
 $NS3_2)$&
 $\smat{2&-1&-1\\
-1&1&-1\\
-3&-1&2}$,& $\smat{
1:\fosp(1|4)\\ 2:G_2\\
3:\fosp(1|4)};
 $&
 \includegraphics{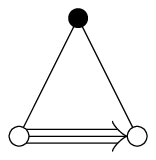}
\end{tabular}

\begin{tabular}
{>{\PBS\raggedright\hspace{0pt}}m{12mm}
 >{\PBS\raggedright\hspace{0pt}}m{27mm}
 >{\PBS\raggedright\hspace{0pt}}m{35mm}
 >{\PBS\raggedright\hspace{0pt}}m{30mm}}
 &&&\\
 $NS3_3)$&
 $\smat{2&-1&-1\\
-1&1&-1\\
-4&-1&2}$,& $\smat{1:\fosp(1|4)\\ 2:A_2^{(2)}\\
3:\fosp(1|4)};
 $&
 \includegraphics{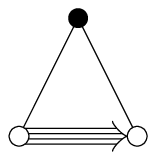}
\end{tabular}

\begin{tabular}
{>{\PBS\raggedright\hspace{0pt}}m{12mm}
 >{\PBS\raggedright\hspace{0pt}}m{27mm}
 >{\PBS\raggedright\hspace{0pt}}m{35mm}
 >{\PBS\raggedright\hspace{0pt}}m{30mm}}
 &&&\\
 $S3_4)$&
 $\smat{2&-1&-2\\
-1&1&-1\\
-2&-1&2}$,& $\smat{1:\fosp(1|4)\\ 2:A_1^{(1)}\\
3:\fosp(1|4)};
 $&
 \includegraphics{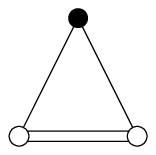}
\end{tabular}

\begin{tabular}
{>{\PBS\raggedright\hspace{0pt}}m{12mm}
 >{\PBS\raggedright\hspace{0pt}}m{27mm}
 >{\PBS\raggedright\hspace{0pt}}m{35mm}
 >{\PBS\raggedright\hspace{0pt}}m{30mm}}
 &&&\\
 $S3_5)$&
 $\smat{2&-1&-1\\ -1&1&-1\\
-1&-1&2}$,& $\smat{1:\fosp(1|4)\\ 2:A_2\\
3:\fosp(1|4)};
 $&
 \includegraphics{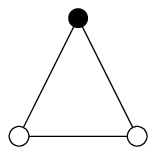}
\end{tabular}

\begin{tabular}
{>{\PBS\raggedright\hspace{0pt}}m{12mm}
 >{\PBS\raggedright\hspace{0pt}}m{27mm}
 >{\PBS\raggedright\hspace{0pt}}m{35mm}
 >{\PBS\raggedright\hspace{0pt}}m{30mm}}
 &&&\\
 $S3_6)$&
 $\smat{2&-1&-1\\
-2&1&-1\\
-2&-1&2}$,& $\smat{1:\fosp(1|4)\\ 2:B_2\\
3:\fosp(1|2)^{(1)}};
 $&
 \includegraphics{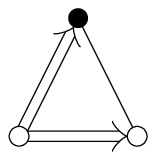}
\end{tabular}

\begin{tabular}
{>{\PBS\raggedright\hspace{0pt}}m{12mm}
 >{\PBS\raggedright\hspace{0pt}}m{27mm}
 >{\PBS\raggedright\hspace{0pt}}m{35mm}
 >{\PBS\raggedright\hspace{0pt}}m{30mm}}
 &&&\\
 $NS3_7)$&
 $\smat{2&-1&-1\\ -2&1&-1\\
-3&-1&2}$,& $\smat{1:\fosp(1|4)\\ 2:G_2\\
3:\fosp(1|2)^{(1)}};
 $&
 \includegraphics{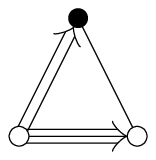}
\end{tabular}

\begin{tabular}
{>{\PBS\raggedright\hspace{0pt}}m{12mm}
 >{\PBS\raggedright\hspace{0pt}}m{27mm}
 >{\PBS\raggedright\hspace{0pt}}m{35mm}
 >{\PBS\raggedright\hspace{0pt}}m{30mm}}
 &&&\\
 $NS3_8)$&
 $\smat{2&-1&-1\\
-2&1&-1\\
-4&-1&2}$,& $\smat{1:\fosp(1|4)\\ 2:A_2^{(2)}\\
3:\fosp(1|2)^{(1)}};$&
 \includegraphics{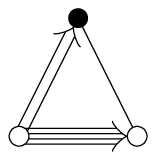}
\end{tabular}

\begin{tabular}
{>{\PBS\raggedright\hspace{0pt}}m{12mm}
 >{\PBS\raggedright\hspace{0pt}}m{27mm}
 >{\PBS\raggedright\hspace{0pt}}m{35mm}
 >{\PBS\raggedright\hspace{0pt}}m{30mm}}
 &&&\\
 $NS3_9)$&
 $\smat{
2&-1&-2\\
-2&1&-1\\
-2&-1&2
}$,& $\smat{1:\fosp(1|4)\\
2:A_1^{(1)}\\
3:\fosp(1|2)^{(1)} };
 $&
 \includegraphics{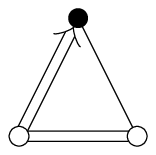}
\end{tabular}

\begin{tabular}
{>{\PBS\raggedright\hspace{0pt}}m{14mm}
 >{\PBS\raggedright\hspace{0pt}}m{27mm}
 >{\PBS\raggedright\hspace{0pt}}m{35mm}
 >{\PBS\raggedright\hspace{0pt}}m{30mm}}
 &&&\\
 $NS3_{10})$&
 $\smat{
2&-1&-1\\
-2&1&-1\\
-1&-1&2
}$,& $\smat{1:\fosp(1|4)\\
2:A_2\\
3:\fosp(1|2)^{(1)} };
 $&
 \includegraphics{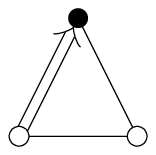}
\end{tabular}

\begin{tabular}
{>{\PBS\raggedright\hspace{0pt}}m{14mm}
 >{\PBS\raggedright\hspace{0pt}}m{27mm}
 >{\PBS\raggedright\hspace{0pt}}m{35mm}
 >{\PBS\raggedright\hspace{0pt}}m{30mm}}
 &&&\\
 $NS3_{11})$&
 $\smat{
2&-1&-2\\
-2&1&-1\\
-1&-1&2
}$,& $\smat{1:\fosp(1|4)\\
2:B_2\\
3:\fosp(1|2)^{(1)} };
 $&
 \includegraphics{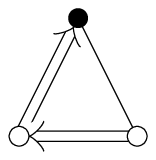}
\end{tabular}

\begin{tabular}
{>{\PBS\raggedright\hspace{0pt}}m{14mm}
 >{\PBS\raggedright\hspace{0pt}}m{27mm}
 >{\PBS\raggedright\hspace{0pt}}m{35mm}
 >{\PBS\raggedright\hspace{0pt}}m{30mm}}
 &&&\\
 $NS3_{12})$&
 $\smat{
2&-1&-3\\
-2&1&-1\\
-1&-1&2}$,& $\smat{
1:\fosp(1|4)\\
2:G_2\\
3:\fosp(1|2)^{(1)} };
 $&
 \includegraphics{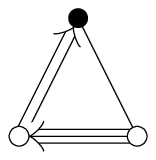}
\end{tabular}

\begin{tabular}
{>{\PBS\raggedright\hspace{0pt}}m{14mm}
 >{\PBS\raggedright\hspace{0pt}}m{27mm}
 >{\PBS\raggedright\hspace{0pt}}m{35mm}
 >{\PBS\raggedright\hspace{0pt}}m{30mm}}
 &&&\\
 $NS3_{13})$&
 $\smat{
2&-1&-4\\
-2&1&-1\\
-1&-1&2
}$,& $\smat{1:\fosp(1|4)\\
2:A_2^{(2)}\\
3:\fosp(1|2)^{(1)}; } $&
 \includegraphics{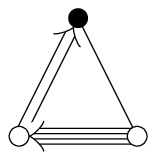}
\end{tabular}

\begin{tabular}
{>{\PBS\raggedright\hspace{0pt}}m{14mm}
 >{\PBS\raggedright\hspace{0pt}}m{27mm}
 >{\PBS\raggedright\hspace{0pt}}m{35mm}
 >{\PBS\raggedright\hspace{0pt}}m{30mm}}
 &&&\\
 $NS3_{14})$&
 $\smat{
2&-2&-1\\
-1&1&-1\\
-2&-1&2
}$,& $\smat{1:\fosp(1|4)\\
2:B_2\\
3:\fsl(1|3)^{(4)} };
 $&
 \includegraphics{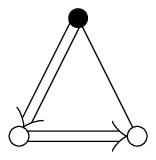}
\end{tabular}

\begin{tabular}
{>{\PBS\raggedright\hspace{0pt}}m{14mm}
 >{\PBS\raggedright\hspace{0pt}}m{27mm}
 >{\PBS\raggedright\hspace{0pt}}m{35mm}
 >{\PBS\raggedright\hspace{0pt}}m{30mm}}
 &&&\\
 $NS3_{15})$&
 $\smat{
2&-2&-1\\
-1&1&-1\\
-3&-1&2
}$,& $\smat{1:\fosp(1|4)\\
2:G_2\\
3:\fsl(1|3)^{(4)} };
 $&
 \includegraphics{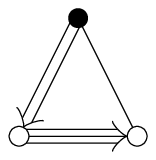}
\end{tabular}

\begin{tabular}
{>{\PBS\raggedright\hspace{0pt}}m{14mm}
 >{\PBS\raggedright\hspace{0pt}}m{27mm}
 >{\PBS\raggedright\hspace{0pt}}m{35mm}
 >{\PBS\raggedright\hspace{0pt}}m{30mm}}
 &&&\\
 $NS3_{16})$&
 $\smat{
2&-2&-1\\
-1&1&-1\\
-4&-1&2
}$,& $\smat{1:\fosp(1|4)\\
2:A_2^{(2)}\\
3:\fsl(1|3)^{(4)} };
 $&
 \includegraphics{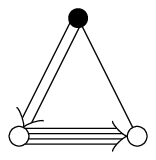}
\end{tabular}

\begin{tabular}
{>{\PBS\raggedright\hspace{0pt}}m{14mm}
 >{\PBS\raggedright\hspace{0pt}}m{27mm}
 >{\PBS\raggedright\hspace{0pt}}m{35mm}
 >{\PBS\raggedright\hspace{0pt}}m{30mm}}
 &&&\\
 $NS3_{17})$&
 $\smat{
2&-2&-2\\
-1&1&-1\\
-2&-1&2
}$,& $\smat{1:\fosp(1|4)\\
2:A_1^{(1)}\\
3:\fsl(1|3)^{(4)} };
 $&
 \includegraphics{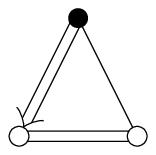}
\end{tabular}

\begin{tabular}
{>{\PBS\raggedright\hspace{0pt}}m{14mm}
 >{\PBS\raggedright\hspace{0pt}}m{27mm}
 >{\PBS\raggedright\hspace{0pt}}m{35mm}
 >{\PBS\raggedright\hspace{0pt}}m{30mm}}
 &&&\\
 $NS3_{18})$&
 $\smat{
2&-2&-1\\
-1&1&-1\\
-1&-1&2
}$,& $\smat{1:\fosp(1|4)\\
2:A_2\\
3:\fsl(1|3)^{(4)} };
 $&
 \includegraphics{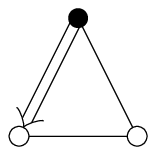}
\end{tabular}

\begin{tabular}
{>{\PBS\raggedright\hspace{0pt}}m{14mm}
 >{\PBS\raggedright\hspace{0pt}}m{27mm}
 >{\PBS\raggedright\hspace{0pt}}m{35mm}
 >{\PBS\raggedright\hspace{0pt}}m{30mm}}
 &&&\\
 $S3_{19})$&
 $\smat{
2&-2&-2\\
-1&1&-1\\
-1&-1&2
}$,& $\smat{1:\fosp(1|4)\\
2:B_2\\
3:\fsl(1|3)^{(4)} };
 $&
 \includegraphics{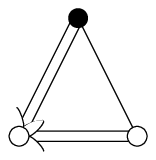}
\end{tabular}

\begin{tabular}
{>{\PBS\raggedright\hspace{0pt}}m{14mm}
 >{\PBS\raggedright\hspace{0pt}}m{27mm}
 >{\PBS\raggedright\hspace{0pt}}m{35mm}
 >{\PBS\raggedright\hspace{0pt}}m{30mm}}
 &&&\\
 $NS3_{20})$&
 $\smat{
2&-2&-3\\
-1&1&-1\\
-1&-1&2
}$,& $\smat{1:\fosp(1|4)\\
2:G_2\\
3:\fsl(1|3)^{(4)} };
 $&
 \includegraphics{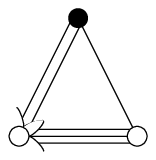}
\end{tabular}

\begin{tabular}
{>{\PBS\raggedright\hspace{0pt}}m{14mm}
 >{\PBS\raggedright\hspace{0pt}}m{27mm}
 >{\PBS\raggedright\hspace{0pt}}m{35mm}
 >{\PBS\raggedright\hspace{0pt}}m{30mm}}
 &&&\\
 $NS3_{21})$&
 $\smat{
2&-2&-4\\
-1&1&-1\\
-1&-1&2
}$,& $\smat{1:\fosp(1|4)\\
2:A_2^{(2)}\\
3:\fsl(1|3)^{(4)} };
 $&
 \includegraphics{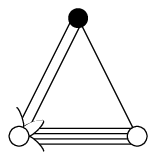}
\end{tabular}

\begin{tabular}
{>{\PBS\raggedright\hspace{0pt}}m{14mm}
 >{\PBS\raggedright\hspace{0pt}}m{27mm}
 >{\PBS\raggedright\hspace{0pt}}m{35mm}
 >{\PBS\raggedright\hspace{0pt}}m{30mm}}
 &&&\\
 $NS3_{22})$&
 $\smat{
2&-1&-1\\
-2&1&-2\\
-2&-1&2
}$,& $\smat{1:\fosp(1|2)^{(1)}\\
2:B_2\\
3:\fosp(1|2)^{(1)} };
 $&
 \includegraphics{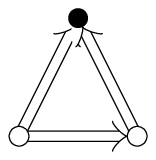}
\end{tabular}

\begin{tabular}
{>{\PBS\raggedright\hspace{0pt}}m{14mm}
 >{\PBS\raggedright\hspace{0pt}}m{27mm}
 >{\PBS\raggedright\hspace{0pt}}m{35mm}
 >{\PBS\raggedright\hspace{0pt}}m{30mm}}
 &&&\\
 $NS3_{23})$&
 $\smat{
2&-1&-1\\
-2&1&-2\\
-3&-1&2
}$,& $\smat{1:\fosp(1|2)^{(1)}\\
2:G_2\\
3:\fosp(1|2)^{(1)} };
 $&
 \includegraphics{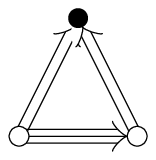}
\end{tabular}

\begin{tabular}
{>{\PBS\raggedright\hspace{0pt}}m{14mm}
 >{\PBS\raggedright\hspace{0pt}}m{27mm}
 >{\PBS\raggedright\hspace{0pt}}m{35mm}
 >{\PBS\raggedright\hspace{0pt}}m{30mm}}
 &&&\\
 $NS3_{24})$&
 $\smat{
2&-1&-1\\
-2&1&-2\\
-4&-1&2
}$,& $\smat{1:\fosp(1|2)^{(1)}\\
2:A_2^{(2)}\\
3:\fosp(1|2)^{(1)} };
 $&
 \includegraphics{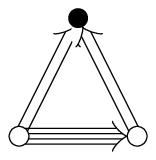}
\end{tabular}

\begin{tabular}
{>{\PBS\raggedright\hspace{0pt}}m{14mm}
 >{\PBS\raggedright\hspace{0pt}}m{27mm}
 >{\PBS\raggedright\hspace{0pt}}m{35mm}
 >{\PBS\raggedright\hspace{0pt}}m{30mm}}
 &&&\\
 $S3_{25})$&
 $\smat{
2&-1&-2\\
-2&1&-2\\
-2&-1&2
}$,& $\smat{1:\fosp(1|2)^{(1)}\\
2:A_1^{(1)}\\
3:\fosp(1|2)^{(1)} };
 $&
 \includegraphics{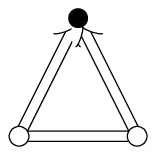}
\end{tabular}

\begin{tabular}
{>{\PBS\raggedright\hspace{0pt}}m{14mm}
 >{\PBS\raggedright\hspace{0pt}}m{27mm}
 >{\PBS\raggedright\hspace{0pt}}m{35mm}
 >{\PBS\raggedright\hspace{0pt}}m{30mm}}
 &&&\\
 $S3_{26})$&
 $\smat{
2&-1&-1\\
-2&1&-2\\
-1&-1&2
}$,& $\smat{1:\fosp(1|2)^{(1)}\\
2:A_2\\
3:\fosp(1|2)^{(1)} };
 $&
 \includegraphics{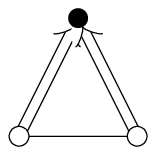}
\end{tabular}

\begin{tabular}
{>{\PBS\raggedright\hspace{0pt}}m{14mm}
 >{\PBS\raggedright\hspace{0pt}}m{27mm}
 >{\PBS\raggedright\hspace{0pt}}m{35mm}
 >{\PBS\raggedright\hspace{0pt}}m{30mm}}
 &&&\\
 $NS3_{27})$&
 $\smat{
2&-2&-1\\
-1&1&-2\\
-2&-1&2
}$,& $\smat{1:\fosp(1|2)^{(1)}\\
2:B_2\\
3:\fsl(1|3)^{(4)} };
 $&
 \includegraphics{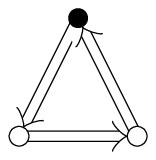}
\end{tabular}

\begin{tabular}
{>{\PBS\raggedright\hspace{0pt}}m{14mm}
 >{\PBS\raggedright\hspace{0pt}}m{27mm}
 >{\PBS\raggedright\hspace{0pt}}m{35mm}
 >{\PBS\raggedright\hspace{0pt}}m{30mm}}
 &&&\\
 $NS3_{28})$&
 $\smat{
2&-2&-1\\
-1&1&-2\\
-3&-1&2
}$,& $\smat{1:\fosp(1|2)^{(1)}\\
2:G_2\\
3:\fsl(1|3)^{(4)} };
 $&
 \includegraphics{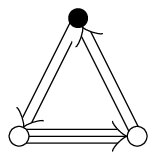}
\end{tabular}

\begin{tabular}
{>{\PBS\raggedright\hspace{0pt}}m{14mm}
 >{\PBS\raggedright\hspace{0pt}}m{27mm}
 >{\PBS\raggedright\hspace{0pt}}m{35mm}
 >{\PBS\raggedright\hspace{0pt}}m{30mm}}
 &&&\\
 $NS3_{29})$&
 $\smat{
2&-2&-1\\
-1&1&-2\\
-4&-1&2
}$,& $\smat{1:\fosp(1|2)^{(1)}\\
2:A_2^{(2)}\\
3:\fsl(1|3)^{(4)} };
 $&
 \includegraphics{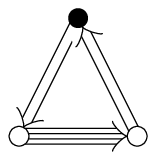}
\end{tabular}

\begin{tabular}
{>{\PBS\raggedright\hspace{0pt}}m{14mm}
 >{\PBS\raggedright\hspace{0pt}}m{27mm}
 >{\PBS\raggedright\hspace{0pt}}m{35mm}
 >{\PBS\raggedright\hspace{0pt}}m{30mm}}
 &&&\\
 $NS3_{30})$&
 $\smat{
2&-2&-2\\
-1&1&-2\\
-2&-1&2
}$,& $\smat{1:\fosp(1|2)^{(1)}\\
2:A_1^{(1)}\\
3:\fsl(1|3)^{(4)} };
 $&
 \includegraphics{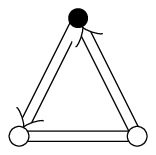}
\end{tabular}

\begin{tabular}
{>{\PBS\raggedright\hspace{0pt}}m{14mm}
 >{\PBS\raggedright\hspace{0pt}}m{27mm}
 >{\PBS\raggedright\hspace{0pt}}m{35mm}
 >{\PBS\raggedright\hspace{0pt}}m{30mm}}
 &&&\\
 $NS3_{31})$&
 $\smat{
2&-2&-1\\
-1&1&-2\\
-1&-1&2
}$,& $\smat{1:\fosp(1|2)^{(1)}\\
2:A_2\\
3:\fsl(1|3)^{(4)} };
 $&
 \includegraphics{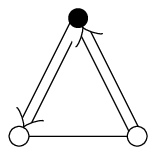}
\end{tabular}

\begin{tabular}
{>{\PBS\raggedright\hspace{0pt}}m{14mm}
 >{\PBS\raggedright\hspace{0pt}}m{27mm}
 >{\PBS\raggedright\hspace{0pt}}m{35mm}
 >{\PBS\raggedright\hspace{0pt}}m{30mm}}
 &&&\\
 $NS3_{32})$&
 $\smat{
2&-2&-2\\
-1&1&-2\\
-1&-1&2
}$,& $\smat{1:\fosp(1|2)^{(1)}\\
2:B_2\\
3:\fsl(1|3)^{(4)} };
 $&
 \includegraphics{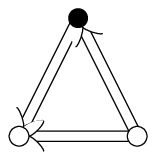}
\end{tabular}

\begin{tabular}
{>{\PBS\raggedright\hspace{0pt}}m{14mm}
 >{\PBS\raggedright\hspace{0pt}}m{27mm}
 >{\PBS\raggedright\hspace{0pt}}m{35mm}
 >{\PBS\raggedright\hspace{0pt}}m{30mm}}
 &&&\\
 $NS3_{33})$&
 $\smat{
2&-2&-3\\
-1&1&-2\\
-1&-1&2
}$,& $\smat{1:\fosp(1|2)^{(1)}\\
2:G_2\\
3:\fsl(1|3)^{(4)} };
 $&
 \includegraphics{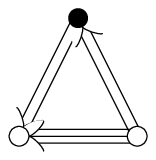}
\end{tabular}

\begin{tabular}
{>{\PBS\raggedright\hspace{0pt}}m{14mm}
 >{\PBS\raggedright\hspace{0pt}}m{27mm}
 >{\PBS\raggedright\hspace{0pt}}m{35mm}
 >{\PBS\raggedright\hspace{0pt}}m{30mm}}
 &&&\\
 $S3_{34})$&
 $\smat{
2&-2&-4\\
-1&1&-2\\
-1&-1&2
}$,& $\smat{1:\fosp(1|2)^{(1)}\\
2:A_2^{(2)}\\
3:\fsl(1|3)^{(4)} };
 $&
 \includegraphics{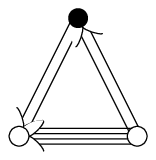}
\end{tabular}

\begin{tabular}
{>{\PBS\raggedright\hspace{0pt}}m{14mm}
 >{\PBS\raggedright\hspace{0pt}}m{27mm}
 >{\PBS\raggedright\hspace{0pt}}m{35mm}
 >{\PBS\raggedright\hspace{0pt}}m{30mm}}
 &&&\\
 $S3_{35})$&
 $\smat{
2&-2&-1\\
-1&1&-1\\
-2&-2&2
}$,& $\smat{1:\fsl(1|3)^{(4)}\\
2:B_2\\
3:\fsl(1|3)^{(4)} };
 $&
 \includegraphics{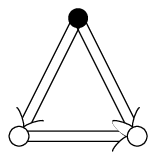}
\end{tabular}

\begin{tabular}
{>{\PBS\raggedright\hspace{0pt}}m{14mm}
 >{\PBS\raggedright\hspace{0pt}}m{27mm}
 >{\PBS\raggedright\hspace{0pt}}m{35mm}
 >{\PBS\raggedright\hspace{0pt}}m{30mm}}
 &&&\\
 $NS3_{36})$&
 $\smat{
2&-2&-1\\
-1&1&-1\\
-3&-2&2
}$,& $\smat{1:\fsl(1|3)^{(4)}\\
2:G_2\\
3:\fsl(1|3)^{(4)} };
 $&
 \includegraphics{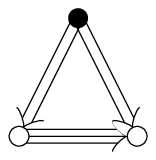}
\end{tabular}

\begin{tabular}
{>{\PBS\raggedright\hspace{0pt}}m{14mm}
 >{\PBS\raggedright\hspace{0pt}}m{27mm}
 >{\PBS\raggedright\hspace{0pt}}m{35mm}
 >{\PBS\raggedright\hspace{0pt}}m{30mm}}
 &&&\\
 $NS3_{37})$&
 $\smat{
2&-2&-1\\
-1&1&-1\\
-4&-2&2
}$,& $\smat{1:\fsl(1|3)^{(4)}\\
2:A_2^{(2)}\\
3:\fsl(1|3)^{(4)} };
 $&
 \includegraphics{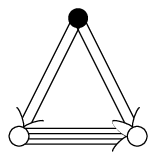}
\end{tabular}

\begin{tabular}
{>{\PBS\raggedright\hspace{0pt}}m{14mm}
 >{\PBS\raggedright\hspace{0pt}}m{27mm}
 >{\PBS\raggedright\hspace{0pt}}m{35mm}
 >{\PBS\raggedright\hspace{0pt}}m{30mm}}
 &&&\\
 $S3_{38})$&
 $\smat{
2&-2&-2\\
-1&1&-1\\
-2&-2&2
}$,& $\smat{1:\fsl(1|3)^{(4)}\\
2:A_1^{(1)}\\
3:\fsl(1|3)^{(4)} };
 $&
 \includegraphics{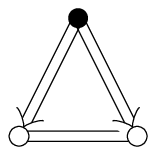}
\end{tabular}

\begin{tabular}
{>{\PBS\raggedright\hspace{0pt}}m{14mm}
 >{\PBS\raggedright\hspace{0pt}}m{27mm}
 >{\PBS\raggedright\hspace{0pt}}m{35mm}
 >{\PBS\raggedright\hspace{0pt}}m{30mm}}
 &&&\\
 $S3_{39})$&
 $\smat{
2&-2&-1\\
-1&1&-1\\
-1&-2&2
}$,& $\smat{1:\fsl(1|3)^{(4)}\\
2:A_2\\
3:\fsl(1|3)^{(4)} };
 $&
 \includegraphics{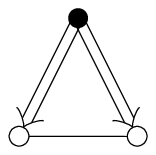}
\end{tabular}

\begin{tabular}
{>{\PBS\raggedright\hspace{0pt}}m{14mm}
 >{\PBS\raggedright\hspace{0pt}}m{27mm}
 >{\PBS\raggedright\hspace{0pt}}m{35mm}
 >{\PBS\raggedright\hspace{0pt}}m{30mm}}
 &&&\\
 $S3_{40})$&
 $\smat{
2&-1&-1\\
-1&1&-1\\
-1&-1&1
}$,& $\smat{1:\fosp(2|2)^{(2)}\\
2:\fosp(1|4)\\
3:\fosp(1|4) };
 $&
 \includegraphics{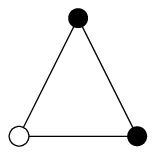}
\end{tabular}

\begin{tabular}
{>{\PBS\raggedright\hspace{0pt}}m{14mm}
 >{\PBS\raggedright\hspace{0pt}}m{27mm}
 >{\PBS\raggedright\hspace{0pt}}m{35mm}
 >{\PBS\raggedright\hspace{0pt}}m{30mm}}
 &&&\\
 $NS3_{41})$&
 $\smat{
2&-1&-1\\
-2&1&-1\\
-1&-1&1
}$,& $\smat{1:\fosp(2|2)^{(2)}\\
2:\fosp(1|4)\\
3:\fosp(1|2)^{(1)} };
 $&
 \includegraphics{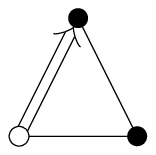}
\end{tabular}

\begin{tabular}
{>{\PBS\raggedright\hspace{0pt}}m{14mm}
 >{\PBS\raggedright\hspace{0pt}}m{27mm}
 >{\PBS\raggedright\hspace{0pt}}m{35mm}
 >{\PBS\raggedright\hspace{0pt}}m{30mm}}
 &&&\\
 $NS3_{42})$&
 $\smat{
2&-2&-1\\
-1&1&-1\\
-1&-1&1
}$,& $\smat{1:\fosp(2|2)^{(2)}\\
2:\fosp(1|4)\\
3:\fsl(1|3)^{(4)} };
 $&
 \includegraphics{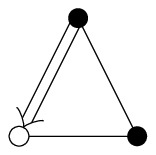}
\end{tabular}

\begin{tabular}
{>{\PBS\raggedright\hspace{0pt}}m{14mm}
 >{\PBS\raggedright\hspace{0pt}}m{27mm}
 >{\PBS\raggedright\hspace{0pt}}m{35mm}
 >{\PBS\raggedright\hspace{0pt}}m{30mm}}
 &&&\\
 $S3_{43})$&
 $\smat{
2&-1&-1\\
-2&1&-1\\
-2&-1&1
}$,& $\smat{1:\fosp(2|2)^{(2)}\\
2:\fosp(1|2)^{(1)}\\
3:\fosp(1|2)^{(1)} };
 $&
 \includegraphics{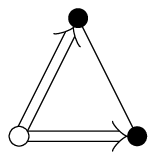}
\end{tabular}

\begin{tabular}
{>{\PBS\raggedright\hspace{0pt}}m{14mm}
 >{\PBS\raggedright\hspace{0pt}}m{27mm}
 >{\PBS\raggedright\hspace{0pt}}m{35mm}
 >{\PBS\raggedright\hspace{0pt}}m{30mm}}
 &&&\\
 $NS3_{44})$&
 $\smat{
2&-2&-1\\
-1&1&-1\\
-2&-1&1
}$,& $\smat{1:\fosp(2|2)^{(2)}\\
2:\fosp(1|2)^{(1)}\\
3:\fsl(1|3)^{(4)} };
 $&
 \includegraphics{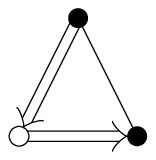}
\end{tabular}

\begin{tabular}
{>{\PBS\raggedright\hspace{0pt}}m{14mm}
 >{\PBS\raggedright\hspace{0pt}}m{27mm}
 >{\PBS\raggedright\hspace{0pt}}m{35mm}
 >{\PBS\raggedright\hspace{0pt}}m{30mm}}
 &&&\\
 $S3_{45})$&
 $\smat{
2&-2&-2\\
-1&1&-1\\
-1&-1&1
}$,& $\smat{1:\fosp(2|2)^{(2)}\\
2:\fsl(1|3)^{(4)}\\
3:\fsl(1|3)^{(4)} };
 $&
 \includegraphics{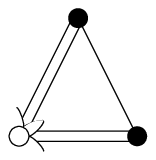}
\end{tabular}

\begin{tabular}
{>{\PBS\raggedright\hspace{0pt}}m{14mm}
 >{\PBS\raggedright\hspace{0pt}}m{27mm}
 >{\PBS\raggedright\hspace{0pt}}m{35mm}
 >{\PBS\raggedright\hspace{0pt}}m{30mm}}
 &&&\\
 $S3_{46})$&
 $\smat{
1&-1&-1\\
-1&1&-1\\
-1&-1&1
}$,& $\smat{1:\fosp(2|2)^{(2)}\\
2:\fosp(2|2)^{(2)}\\
3:\fosp(2|2)^{(2)} };
 $&
 \includegraphics{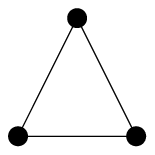}
\end{tabular}

\begin{tabular}
{>{\PBS\raggedright\hspace{0pt}}m{14mm}
 >{\PBS\raggedright\hspace{0pt}}m{27mm}
 >{\PBS\raggedright\hspace{0pt}}m{35mm}
 >{\PBS\raggedright\hspace{0pt}}m{30mm}}
 &&&\\
 $S3_{47})$&
 $\smat{
2&-1&0\\
-2&1&-1\\
0&-1&2
}$,& $\smat{1:\fosp(1|4)\\
2:\\
3:\fosp(1|2)^{(1)} };
 $&
 \includegraphics{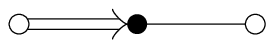}
\end{tabular}

\begin{tabular}
{>{\PBS\raggedright\hspace{0pt}}m{14mm}
 >{\PBS\raggedright\hspace{0pt}}m{27mm}
 >{\PBS\raggedright\hspace{0pt}}m{35mm}
 >{\PBS\raggedright\hspace{0pt}}m{30mm}}
 &&&\\
 $S3_{48})$&
 $\smat{
2&-2&0\\
-1&1&-1\\
0&-1&2
}$,& $\smat{1:\fosp(1|4)\\
2:\\
3:\fsl(1|3)^{(4)} };
 $&
 \includegraphics{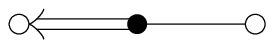}
\end{tabular}

\begin{tabular}
{>{\PBS\raggedright\hspace{0pt}}m{14mm}
 >{\PBS\raggedright\hspace{0pt}}m{27mm}
 >{\PBS\raggedright\hspace{0pt}}m{35mm}
 >{\PBS\raggedright\hspace{0pt}}m{30mm}}
 &&&\\
 $S3_{49})$&
 $\smat{
2&-1&0\\
-2&1&-2\\
0&-1&2
}$,& $\smat{1:\fosp(1|2)^{(1)}\\
2:\\
3:\fosp(1|2)^{(1)} };
 $&
 \includegraphics{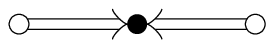}
\end{tabular}

\begin{tabular}
{>{\PBS\raggedright\hspace{0pt}}m{14mm}
 >{\PBS\raggedright\hspace{0pt}}m{27mm}
 >{\PBS\raggedright\hspace{0pt}}m{35mm}
 >{\PBS\raggedright\hspace{0pt}}m{30mm}}
 &&&\\
 $S3_{50})$&
 $\smat{
2&-2&0\\
-1&1&-2\\
0&-1&2
}$,& $\smat{1:\fosp(1|2)^{(1)}\\
2:\\
3:\fsl(1|3)^{(4)} };
 $&
 \includegraphics{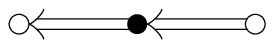}
\end{tabular}

\begin{tabular}
{>{\PBS\raggedright\hspace{0pt}}m{14mm}
 >{\PBS\raggedright\hspace{0pt}}m{27mm}
 >{\PBS\raggedright\hspace{0pt}}m{35mm}
 >{\PBS\raggedright\hspace{0pt}}m{30mm}}
 &&&\\
 $S3_{51})$&
 $\smat{
2&-2&0\\
-1&1&-1\\
0&-2&2
}$,& $\smat{1:\fsl(1|3)^{(4)}\\
2:\\
3:\fsl(1|3)^{(4)} };
 $&
 \includegraphics{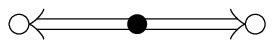}
\end{tabular}

\begin{tabular}
{>{\PBS\raggedright\hspace{0pt}}m{14mm}
 >{\PBS\raggedright\hspace{0pt}}m{27mm}
 >{\PBS\raggedright\hspace{0pt}}m{35mm}
 >{\PBS\raggedright\hspace{0pt}}m{30mm}}
 &&&\\
 $S3_{52})$&
 $\smat{
2&-1&0\\
-3&2&-1\\
0&-1&1
}$,& $\smat{1:\fosp(1|4)\\
2:\\
3:G_2 };
 $&
 \includegraphics{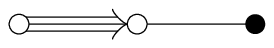}
\end{tabular}

\begin{tabular}
{>{\PBS\raggedright\hspace{0pt}}m{14mm}
 >{\PBS\raggedright\hspace{0pt}}m{27mm}
 >{\PBS\raggedright\hspace{0pt}}m{35mm}
 >{\PBS\raggedright\hspace{0pt}}m{30mm}}
 &&&\\
 $S3_{53})$&
 $\smat{
2&-1&0\\
-4&2&-1\\
0&-1&1
}$,& $\smat{1:\fosp(1|4)\\
2:\\
3:A_2^{(2)} };
 $&
 \includegraphics{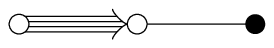}
\end{tabular}

\begin{tabular}
{>{\PBS\raggedright\hspace{0pt}}m{14mm}
 >{\PBS\raggedright\hspace{0pt}}m{27mm}
 >{\PBS\raggedright\hspace{0pt}}m{35mm}
 >{\PBS\raggedright\hspace{0pt}}m{30mm}}
 &&&\\
 $S3_{54})$&
 $\smat{
2&-2&0\\
-2&2&-1\\
0&-1&1
}$,& $\smat{1:\fosp(1|4)\\
2:\\
3:A_1^{(1)} };
 $&
 \includegraphics{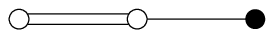}
\end{tabular}

\begin{tabular}
{>{\PBS\raggedright\hspace{0pt}}m{14mm}
 >{\PBS\raggedright\hspace{0pt}}m{27mm}
 >{\PBS\raggedright\hspace{0pt}}m{35mm}
 >{\PBS\raggedright\hspace{0pt}}m{30mm}}
 &&&\\
 $S3_{55})$&
 $\smat{
2&-3&0\\
-1&2&-1\\
0&-1&1
}$,& $\smat{1:\fosp(1|4)\\
2:\\
3:G_2 };
 $&
 \includegraphics{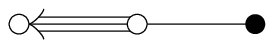}
\end{tabular}

\begin{tabular}
{>{\PBS\raggedright\hspace{0pt}}m{14mm}
 >{\PBS\raggedright\hspace{0pt}}m{27mm}
 >{\PBS\raggedright\hspace{0pt}}m{35mm}
 >{\PBS\raggedright\hspace{0pt}}m{30mm}}
 &&&\\
 $S3_{56})$&
 $\smat{
2&-4&0\\
-1&2&-1\\
0&-1&1
}$,& $\smat{1:\fosp(1|4)\\
2:\\
3:A_2^{(2)} };
 $&
 \includegraphics{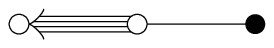}
\end{tabular}

\begin{tabular}
{>{\PBS\raggedright\hspace{0pt}}m{14mm}
 >{\PBS\raggedright\hspace{0pt}}m{27mm}
 >{\PBS\raggedright\hspace{0pt}}m{35mm}
 >{\PBS\raggedright\hspace{0pt}}m{30mm}}
 &&&\\
 $S3_{57})$&
 $\smat{
2&-1&0\\
-2&2&-1\\
0&-2&1
}$,& $\smat{1:\fosp(1|2)^{(1)}\\
2:\\
3:B_2 };
 $&
 \includegraphics{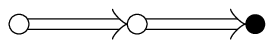}
\end{tabular}

\begin{tabular}
{>{\PBS\raggedright\hspace{0pt}}m{14mm}
 >{\PBS\raggedright\hspace{0pt}}m{27mm}
 >{\PBS\raggedright\hspace{0pt}}m{35mm}
 >{\PBS\raggedright\hspace{0pt}}m{30mm}}
 &&&\\
 $S3_{58})$&
 $\smat{
2&-1&0\\
-3&2&-1\\
0&-2&1
}$,& $\smat{1:\fosp(1|2)^{(1)}\\
2:\\
3:G_2 };
 $&
 \includegraphics{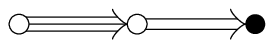}
\end{tabular}

\begin{tabular}
{>{\PBS\raggedright\hspace{0pt}}m{14mm}
 >{\PBS\raggedright\hspace{0pt}}m{27mm}
 >{\PBS\raggedright\hspace{0pt}}m{35mm}
 >{\PBS\raggedright\hspace{0pt}}m{30mm}}
 &&&\\
 $S3_{59})$&
 $\smat{
2&-1&0\\
-4&2&-1\\
0&-2&1
}$,& $\smat{1:\fosp(1|2)^{(1)}\\
2:\\
3:A_2^{(2)} };
 $&
 \includegraphics{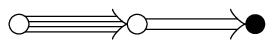}
\end{tabular}

\begin{tabular}
{>{\PBS\raggedright\hspace{0pt}}m{14mm}
 >{\PBS\raggedright\hspace{0pt}}m{27mm}
 >{\PBS\raggedright\hspace{0pt}}m{35mm}
 >{\PBS\raggedright\hspace{0pt}}m{30mm}}
 &&&\\
 $S3_{60})$&
 $\smat{
2&-2&0\\
-2&2&-1\\
0&-2&1
}$,& $\smat{1:\fosp(1|2)^{(1)}\\
2:\\
3:A_1^{(1)} };
 $&
 \includegraphics{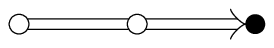}
\end{tabular}

\begin{tabular}
{>{\PBS\raggedright\hspace{0pt}}m{14mm}
 >{\PBS\raggedright\hspace{0pt}}m{27mm}
 >{\PBS\raggedright\hspace{0pt}}m{35mm}
 >{\PBS\raggedright\hspace{0pt}}m{30mm}}
 &&&\\
 $S3_{61})$&
 $\smat{
2&-1&0\\
-1&2&-1\\
0&-2&1
}$,& $\smat{1:\fosp(1|2)^{(1)}\\
2:\\
3:A_2 };
 $&
 \includegraphics{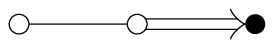}
\end{tabular}

\begin{tabular}
{>{\PBS\raggedright\hspace{0pt}}m{14mm}
 >{\PBS\raggedright\hspace{0pt}}m{27mm}
 >{\PBS\raggedright\hspace{0pt}}m{35mm}
 >{\PBS\raggedright\hspace{0pt}}m{30mm}}
 &&&\\
 $S3_{62})$&
 $\smat{
2&-2&0\\
-1&2&-1\\
0&-2&1
}$,& $\smat{1:\fosp(1|2)^{(1)}\\
2:\\
3:B_2 };
 $&
 \includegraphics{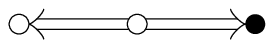}
\end{tabular}

\begin{tabular}
{>{\PBS\raggedright\hspace{0pt}}m{14mm}
 >{\PBS\raggedright\hspace{0pt}}m{27mm}
 >{\PBS\raggedright\hspace{0pt}}m{35mm}
 >{\PBS\raggedright\hspace{0pt}}m{30mm}}
 &&&\\
 $S3_{63})$&
 $\smat{
2&-3&0\\
-1&2&-1\\
0&-2&1
}$,& $\smat{1:\fosp(1|2)^{(1)}\\
2:\\
3:G_2 };
 $&
 \includegraphics{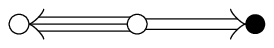}
\end{tabular}

\begin{tabular}
{>{\PBS\raggedright\hspace{0pt}}m{14mm}
 >{\PBS\raggedright\hspace{0pt}}m{27mm}
 >{\PBS\raggedright\hspace{0pt}}m{35mm}
 >{\PBS\raggedright\hspace{0pt}}m{30mm}}
 &&&\\
 $S3_{64})$&
 $\smat{
2&-4&0\\
-1&2&-1\\
0&-2&1
}$,& $\smat{1:\fosp(1|2)^{(1)}\\
2:\\
3:A_2^{(2)} };
 $&
 \includegraphics{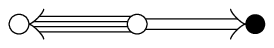}
\end{tabular}

\begin{tabular}
{>{\PBS\raggedright\hspace{0pt}}m{14mm}
 >{\PBS\raggedright\hspace{0pt}}m{27mm}
 >{\PBS\raggedright\hspace{0pt}}m{35mm}
 >{\PBS\raggedright\hspace{0pt}}m{30mm}}
 &&&\\
 $S3_{65})$&
 $\smat{
2&-1&0\\
-2&2&-2\\
0&-1&1
}$,& $\smat{1:\fsl(1|3)^{(4)}\\
2:\\
3:B_2 };
 $&
 \includegraphics{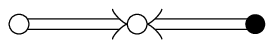}
\end{tabular}

\begin{tabular}
{>{\PBS\raggedright\hspace{0pt}}m{14mm}
 >{\PBS\raggedright\hspace{0pt}}m{27mm}
 >{\PBS\raggedright\hspace{0pt}}m{35mm}
 >{\PBS\raggedright\hspace{0pt}}m{30mm}}
 &&&\\
 $S3_{66})$&
 $\smat{
2&-1&0\\
-3&2&-2\\
0&-1&1
}$,& $\smat{1:\fsl(1|3)^{(4)}\\
2:\\
3:G_2 };
 $&
 \includegraphics{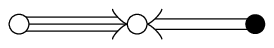}
\end{tabular}

\begin{tabular}
{>{\PBS\raggedright\hspace{0pt}}m{14mm}
 >{\PBS\raggedright\hspace{0pt}}m{27mm}
 >{\PBS\raggedright\hspace{0pt}}m{35mm}
 >{\PBS\raggedright\hspace{0pt}}m{30mm}}
 &&&\\
 $S3_{67})$&
 $\smat{
2&-1&0\\
-4&2&-2\\
0&-1&1
}$,& $\smat{1:\fsl(1|3)^{(4)}\\
2:\\
3:A_2^{(2)} };
 $&
 \includegraphics{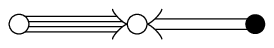}
\end{tabular}

\begin{tabular}
{>{\PBS\raggedright\hspace{0pt}}m{14mm}
 >{\PBS\raggedright\hspace{0pt}}m{27mm}
 >{\PBS\raggedright\hspace{0pt}}m{35mm}
 >{\PBS\raggedright\hspace{0pt}}m{30mm}}
 &&&\\
 $S3_{68})$&
 $\smat{
2&-2&0\\
-2&2&-2\\
0&-1&1
}$,& $\smat{1:\fsl(1|3)^{(4)}\\
2:\\
3:A_1^{(1)} };
 $&
 \includegraphics{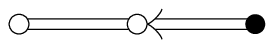}
\end{tabular}

\begin{tabular}
{>{\PBS\raggedright\hspace{0pt}}m{14mm}
 >{\PBS\raggedright\hspace{0pt}}m{27mm}
 >{\PBS\raggedright\hspace{0pt}}m{35mm}
 >{\PBS\raggedright\hspace{0pt}}m{30mm}}
 &&&\\
 $S3_{69})$&
 $\smat{
2&-1&0\\
-1&2&-2\\
0&-1&1
}$,& $\smat{1:\fsl(1|3)^{(4)}\\
2:\\
3:A_2 };
 $&
 \includegraphics{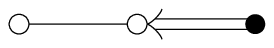}
\end{tabular}

\begin{tabular}
{>{\PBS\raggedright\hspace{0pt}}m{14mm}
 >{\PBS\raggedright\hspace{0pt}}m{27mm}
 >{\PBS\raggedright\hspace{0pt}}m{35mm}
 >{\PBS\raggedright\hspace{0pt}}m{30mm}}
 &&&\\
 $S3_{70})$&
 $\smat{
2&-2&0\\
-1&2&-2\\
0&-1&1
}$,& $\smat{1:\fsl(1|3)^{(4)}\\
2:\\
3:B_2 };
 $&
 \includegraphics{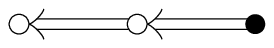}
\end{tabular}

\begin{tabular}
{>{\PBS\raggedright\hspace{0pt}}m{14mm}
 >{\PBS\raggedright\hspace{0pt}}m{27mm}
 >{\PBS\raggedright\hspace{0pt}}m{35mm}
 >{\PBS\raggedright\hspace{0pt}}m{30mm}}
 &&&\\
 $S3_{71})$&
 $\smat{
2&-3&0\\
-1&2&-2\\
0&-1&1
}$,& $\smat{1:\fsl(1|3)^{(4)}\\
2:\\
3:G_2 };
 $&
 \includegraphics{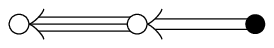}
\end{tabular}

\begin{tabular}
{>{\PBS\raggedright\hspace{0pt}}m{14mm}
 >{\PBS\raggedright\hspace{0pt}}m{27mm}
 >{\PBS\raggedright\hspace{0pt}}m{35mm}
 >{\PBS\raggedright\hspace{0pt}}m{30mm}}
 &&&\\
 $S3_{72})$&
 $\smat{
2&-4&0\\
-1&2&-2\\
0&-1&1
}$,& $\smat{1:\fsl(1|3)^{(4)}\\
2:\\
3:A_2^{(2)} };
 $&
 \includegraphics{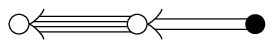}
\end{tabular}

\begin{tabular}
{>{\PBS\raggedright\hspace{0pt}}m{14mm}
 >{\PBS\raggedright\hspace{0pt}}m{27mm}
 >{\PBS\raggedright\hspace{0pt}}m{35mm}
 >{\PBS\raggedright\hspace{0pt}}m{30mm}}
 &&&\\
 $S3_{73})$&
 $\smat{
1&-2&0\\
-1&2&-1\\
0&-1&1
}$,& $\smat{1:\fosp(1|4)\\
2:\\
3:\fosp(1|2)^{(1)} };
 $&
 \includegraphics{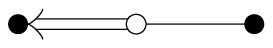}
\end{tabular}

\begin{tabular}
{>{\PBS\raggedright\hspace{0pt}}m{14mm}
 >{\PBS\raggedright\hspace{0pt}}m{27mm}
 >{\PBS\raggedright\hspace{0pt}}m{35mm}
 >{\PBS\raggedright\hspace{0pt}}m{30mm}}
 &&&\\
 $S3_{74})$&
 $\smat{
1&-1&0\\
-2&2&-1\\
0&-1&1
}$,& $\smat{1:\fosp(1|4)\\
2:\\
3:\fsl(1|3)^{(4)} };
 $&
 \includegraphics{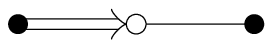}
\end{tabular}

\begin{tabular}
{>{\PBS\raggedright\hspace{0pt}}m{14mm}
 >{\PBS\raggedright\hspace{0pt}}m{27mm}
 >{\PBS\raggedright\hspace{0pt}}m{35mm}
 >{\PBS\raggedright\hspace{0pt}}m{30mm}}
 &&&\\
 $S3_{75})$&
 $\smat{
1&-2&0\\
-1&2&-1\\
0&-2&1
}$,& $\smat{1:\fosp(1|2)^{(1)}\\
2:\\
3:\fosp(1|2)^{(1)} };
 $&
 \includegraphics{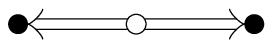}
\end{tabular}

\begin{tabular}
{>{\PBS\raggedright\hspace{0pt}}m{14mm}
 >{\PBS\raggedright\hspace{0pt}}m{27mm}
 >{\PBS\raggedright\hspace{0pt}}m{35mm}
 >{\PBS\raggedright\hspace{0pt}}m{30mm}}
 &&&\\
 $S3_{76})$&
 $\smat{
1&-1&0\\
-2&2&-1\\
0&-2&1
}$,& $\smat{1:\fosp(1|2)^{(1)}\\
2:\\
3:\fsl(1|3)^{(4)} };
 $&
 \includegraphics{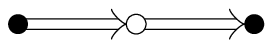}
\end{tabular}

\begin{tabular}
{>{\PBS\raggedright\hspace{0pt}}m{14mm}
 >{\PBS\raggedright\hspace{0pt}}m{27mm}
 >{\PBS\raggedright\hspace{0pt}}m{35mm}
 >{\PBS\raggedright\hspace{0pt}}m{30mm}}
 &&&\\
 $S3_{77})$&
 $\smat{
1&-1&0\\
-2&2&-2\\
0&-1&1
}$,& $\smat{1:\fsl(1|3)^{(4)}\\
2:\\
3:\fsl(1|3)^{(4)} };
 $&
 \includegraphics{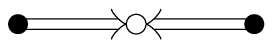}
\end{tabular}

\begin{tabular}
{>{\PBS\raggedright\hspace{0pt}}m{14mm}
 >{\PBS\raggedright\hspace{0pt}}m{27mm}
 >{\PBS\raggedright\hspace{0pt}}m{35mm}
 >{\PBS\raggedright\hspace{0pt}}m{30mm}}
 &&&\\
 $S3_{78})$&
 $\smat{
2&-1&0\\
-1&1&-1\\
0&-1&1
}$,& $\smat{1:\fosp(2|2)^{(2)}\\
2:\\
3:\fosp(1|4) };
 $&
 \includegraphics{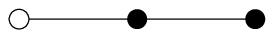}
\end{tabular}

\begin{tabular}
{>{\PBS\raggedright\hspace{0pt}}m{14mm}
 >{\PBS\raggedright\hspace{0pt}}m{27mm}
 >{\PBS\raggedright\hspace{0pt}}m{35mm}
 >{\PBS\raggedright\hspace{0pt}}m{30mm}}
 &&&\\
 $S3_{79})$&
 $\smat{
2&-1&0\\
-2&1&-1\\
0&-1&1
}$,& $\smat{1:\fosp(2|2)^{(2)}\\
2:\\
3:\fosp(1|2)^{(1)} };
 $&
 \includegraphics{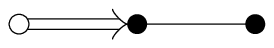}
\end{tabular}

\begin{tabular}
{>{\PBS\raggedright\hspace{0pt}}m{14mm}
 >{\PBS\raggedright\hspace{0pt}}m{27mm}
 >{\PBS\raggedright\hspace{0pt}}m{35mm}
 >{\PBS\raggedright\hspace{0pt}}m{30mm}}
 &&&\\
 $S3_{80})$&
 $\smat{
2&-2&0\\
-1&1&-1\\
0&-1&1
}$,& $\smat{1:\fosp(2|2)^{(2)}\\
2:\\
3:\fsl(1|3)^{(4)} };
 $&
 \includegraphics{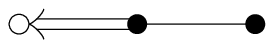}
\end{tabular}

\begin{tabular}
{>{\PBS\raggedright\hspace{0pt}}m{14mm}
 >{\PBS\raggedright\hspace{0pt}}m{27mm}
 >{\PBS\raggedright\hspace{0pt}}m{35mm}
 >{\PBS\raggedright\hspace{0pt}}m{30mm}}
 &&&\\
 $S3_{81})$&
 $\smat{
1&-1&0\\
-1&1&-1\\
0&-1&1
}$,& $\smat{1:\fosp(2|2)^{(2)}\\
2:\\
3:\fosp(2|2)^{(2)} };
 $&
 \includegraphics{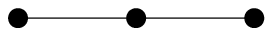}
\end{tabular}

\sssec{Notation: On matrices with a \lq\lq --" sign}\label{recmat}
The rectangular matrix at the beginning of each list of inequivalent
Cartan matrices for each Lie superalgebra shows the result of odd
reflections\footnote{The reader might be interested in analogs of
these results for Lie superalgebras over fields of positive
characteristic, see \cite{BGL}.} (the number of the row is the
number of the Cartan matrix in the list below, the number of the
column is the number of the root (given by small boxed number at the
vertex of the corresponding Dynkin diagram) in which the reflection
is made; the cells contain the results of reflections (the number of
the Cartan matrix obtained) or a \lq\lq --" if the reflection is not
appropriate because $A_{ii}\neq 0$.

For each of the equations given, indicated is one of the two roots.
With each root $a$ its inverse $a^{-1}$ is also a root and to it the
same Lie superalgebra corresponds, its Cartan matrices are obtained
from the given ones by an automorphism or renumbering of the given
matrices.

$NS3_{82})$\quad $5a^2+11a+5=0$,  $a=\frac{\sqrt{21}-11}{10}$
\[


 \includegraphics{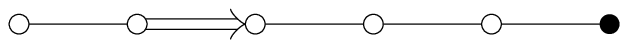}

\begin{tabular}
{>{\PBS\raggedright\hspace{0pt}}m{8mm}
 >{\PBS\raggedright\hspace{0pt}}m{45mm}
 >{\PBS\raggedright\hspace{0pt}}m{42mm}
}
 &&\\
 $S6_{7})$&
 $\smat{
2&-1&0&0&0&0\\
-1&2&-2&0&0&0\\
0&-1&2&-1&0&0\\
0&0&-1&2&-1&0\\
0&0&0&-1&2&-1\\
0&0&0&0&-1&1
}$,& $\smat{1:\fsl(1|9)^{(4)}\\
2:\fosp(1|8)\\
3:\fosp(1|6),A_2\\
4:\fosp(1|4),C_3\\
5:F_4\\
6:E_6^{(2)}};
 $\end{tabular}

 \includegraphics{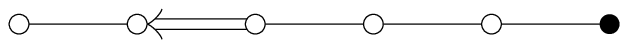}

\ssec{Size 7}{}~{}

\begin{tabular}
{>{\PBS\raggedright\hspace{0pt}}m{8mm}
 >{\PBS\raggedright\hspace{0pt}}m{50mm}
 >{\PBS\raggedright\hspace{0pt}}m{42mm}
}
 &&\\
 $S7_1)$&
 $\smat{2&-1&0&0&0&0&0\\
-1&2&-1&0&0&0&0\\
0&-1&2&-1&-1&0&0\\
0&0&-1&2&0&0&0\\
0&0&-1&0&2&-1&0\\
0&0&0&0&-1&2&-1\\
0&0&0&0&0&-1&1
}$&$\smat{1:\fsl(1|10)^{(2)}\\
2:\fosp(1|10)\\
3:\fosp(1|6),A_2\\
4:\fosp(1|12)\\
5:\fosp(1|4),A_4\\
6:D_5\\
7:E_6};
 $\end{tabular}

 \includegraphics{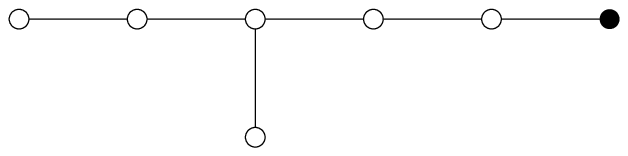}

\ssec{Size 8}{}~{}

\begin{tabular}
{>{\PBS\raggedright\hspace{0pt}}m{8mm}
 >{\PBS\raggedright\hspace{0pt}}m{60mm}
 >{\PBS\raggedright\hspace{0pt}}m{50mm}}
 &&\\
 $S8_1)$&
 $\smat{2&-1&0&0&0&0&0&0\\
-1&2&-1&0&0&0&0&0\\
0&-1&2&-1&-1&0&0&0\\
0&0&-1&2&0&0&0&0\\
0&0&-1&0&2&-1&0&0\\
0&0&0&0&-1&2&-1&0\\
0&0&0&0&0&-1&2&-1\\
0&0&0&0&0&0&-1&1
}$, &$\smat{1:\fsl(1|12)^{(2)}\\
2:\fosp(1|12)\\
3:\fosp(1|8),A_2\\
4:\fosp(1|14)\\
5:\fosp(1|6),A_4\\
6:D_5,\fosp(1|4)\\
7:E_6\\
8:E_7}$
\end{tabular}

\includegraphics{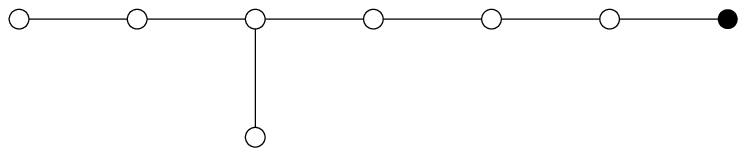}

\ssec{Size 9}{}~{}

\begin{tabular}
{>{\PBS\raggedright\hspace{0pt}}m{8mm}
 >{\PBS\raggedright\hspace{0pt}}m{70mm}
 >{\PBS\raggedright\hspace{0pt}}m{40mm}}
 &&\\
 $S9_1)$&
 $\smat{2&-1&0&0&0&0&0&0&0\\
-1&2&-1&0&0&0&0&0&0\\
0&-1&2&-1&-1&0&0&0&0\\
0&0&-1&2&0&0&0&0&0\\
0&0&-1&0&2&-1&0&0&0\\
0&0&0&0&-1&2&-1&0&0\\
0&0&0&0&0&-1&2&-1&0\\
0&0&0&0&0&0&-1&2&-1\\
0&0&0&0&0&0&0&-1&1
}$&$\smat{1:\fsl(1|14)^{(2)}\\
2:\fosp(1|14)\\
3:\fosp(1|10),A_2\\
4:\fosp(1|16)\\
5:\fosp(1|8),A_4\\
6:D_5,\fosp(1|6)\\
7:E_6,\fosp(1|4)\\
8:E_7\\
9:E_8 }$\end{tabular}%

\ssbegin{Theorem}\label{theor} The above list of $138$ Lie
superalgebras in this section exhausts all almost affine Lie
superalgebras (which are not Lie algebras) of rank
$>2$.\end{Theorem}

In particular, there are no almost affine Lie superalgebras with
either non-symmetrizable Cartan matrices or matrices of larger size.

Observe that to get the Cartan matrix of the even part of any of the
above Lie superalgebras whose Cartan matrix has no 0 on the diagonal
is easy, see \cite{vdL}: just multiply the row with a 1 on the
diagonal by 2. In this way we get 133 almost affine Lie algebras.

\section{Several precise spots where mistakes crept into [FS]}
Since this paper is the third approximation to the true list of
almost affine Lie superalgebras with indecomposable Cartan matrices,
let us point at several of the wrong places in the previous claims.
First, the restriction on the size of the Cartan matrix was wrong
both in \cite{TDP} and \cite{FS}.

The last diagram on p.8 of the arXiv version of \cite{FS} is
manifestly not \lq\lq hyperbolic": Delete the second or third vertex
and check. Same applies to the penultimate diagram; delete the
second vertex. In one the remaining submatrices, in the row with a 0
on the diagonal all entries are of the same sign. This {\it could}
only happen for few values of parameter for $\fosp(4|2;a)$ and
$\fsvect(1|2;a)$ but in reality {\bf may not} even in this case.
Moreover, a bit further, Frappat and Sciarrino list this submatrix
as \lq\lq hyperbolic" ($\#54$ of rank 4) which is an additional
argument against this example.

None of the Cartan matrices with 0 on the diagonal given in
\cite{TDP} and \cite{FS} correspond to any of almost affine Lie
superalgebras.

\section{Appendix: Almost affine Lie algebras} Since the result of
Kobayashi and Morita \cite{KoMo}, and that of Li Wang Lai \cite{Li},
are not easy to get, and Sa{\c{c}}lio{\u{g}}lu's results are
incomplete, we reproduce Li's results in the user-friendly form of
Sa{\c{c}}lio{\u{g}}lu adopted above.




 \includegraphics{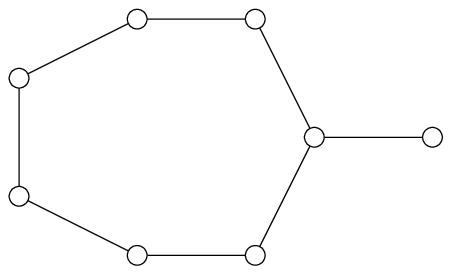}

%

 \includegraphics{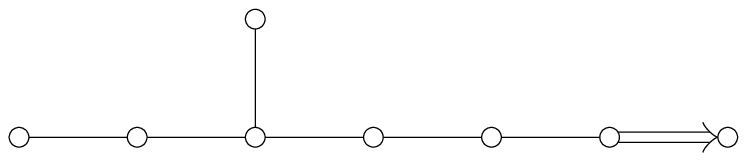}

\newpage

%

 \includegraphics{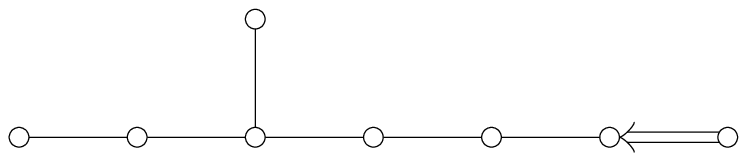}

%

 \includegraphics{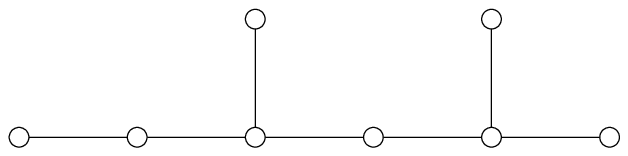}

%

 \includegraphics{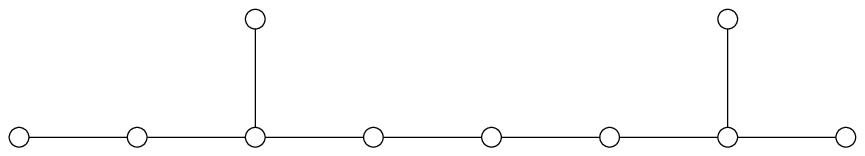}


\end{document}